\documentclass[graybox,natbib]{svmult}
\usepackage{natbib}

\usepackage{type1cm}        
\usepackage{makeidx}         
\usepackage{graphicx}        
\usepackage{multicol}        
\usepackage[bottom]{footmisc}

\usepackage{newtxtext}       %
\usepackage{newtxmath}       

\usepackage{paralist}

\usepackage{pgfplots}
\usepackage{color}
\usepackage{tikz}
\usepackage{pgfplots}
\pgfplotsset{width=8.5cm,compat=1.9}
\usetikzlibrary{intersections,shapes.arrows}

\usepackage{pgfplots}
\usetikzlibrary{matrix}
\usepgfplotslibrary{groupplots}
\pgfplotsset{compat=newest}
\definecolor{grey}{rgb}{0.43, 0.5, 0.5}
\usepackage{amsbsy}
\usepackage{amsfonts}
\usepackage{amsmath}
\providecommand{\bo}{\mathbf}

\providecommand{\normal}{\mathcal{N}}

\renewcommand{\E}{\mathrm{\bo E}}

\newcommand{\n}{^{(n)}}
\newcommand{\pr}{^{\prime}}

\newcommand{\zb}{{\bf z}}
\newcommand{\Zb}{{\bf Z}}
\newcommand{\ub}{{\bf u}}

\newcommand{\Fb}{{\bf F}}

\newcommand{\Qb}{{\bf Q}}
\newcommand{\Sb}{{\bf S}}
\newcommand{\Jb}{{\bf J}}
\newcommand{\Sigmab}{{\boldsymbol\Sigma}}
\newcommand{\Deltab}{{\boldsymbol\Delta}}   

\providecommand{\argmin}{\operatornamewithlimits{\arg\min}}

\usepackage{stackengine}
\stackMath
\newcommand\tenq[2][1]{%
\def\useanchorwidth{T}%
\ifnum#1>1%
\stackunder[0pt]{\tenq[\numexpr#1-1\relax]{#2}}{\scriptscriptstyle\thicksim}%
\else%
\stackunder[1pt]{#2}{\scriptscriptstyle\thicksim}%
\fi%
}

\definecolor{darkraspberry}{rgb}{0.53, 0.15, 0.34}
\definecolor{aurometalsaurus}{rgb}{0.43, 0.5, 0.5}
\definecolor{britishracinggreen}{rgb}{0.0, 0.26, 0.15}
\definecolor{ocre}{RGB}{243,102,25} 
\definecolor{royalblue}{RGB}{0,78,156}
\definecolor{brightgreen}{rgb}{0.4, 1.0, 0.0}

\begin{document}

\title*{On the Finite-Sample Performance of\\  Measure Transportation-Based \\ Multivariate Rank Tests}
\titlerunning{Multivariate rank tests} 
\author{Marc Hallin and Gilles Mordant}
\institute{Marc Hallin \at ECARES and D\'epartement de Math\'ematique, Universit\'e libre de Bruxelles, 
Avenue F.D.\ Roosevelt~50, 1050
 Brussels, Belgium\\  \email{mhallin@ulb.ac.be}
\and Gilles Mordant \at IMS, Universit\"at G\"ottingen, 
 Goldschmidtstra\ss e 7, 37077
 G\"ottingen, Germany\\ \email{mordantgilles@gmail.com}}
%
%
\maketitle

\abstract{Extending to dimension 2 and higher the dual univariate concepts of ranks and quantiles has remained an open problem for more than half a century. Based on measure transportation results, a  solution has been proposed recently under the name {\it center-outward ranks and quantiles} which, contrary to previous proposals, enjoys all the properties that make univariate ranks a successful tool for statistical inference. Just as their univariate counterparts (to which they reduce in dimension one), center-outward ranks allow for the construction of distribution-free and asymptotically efficient tests for a variety of problems where the density of some noise or innovation remains unspecified. The actual implementation of these tests involves the somewhat arbitrary choice of a grid. While the asymptotic impact of that choice is nil, its  finite-sample consequences are not. In this note, we investigate the  finite-sample impact of that choice in the typical context of the multivariate two-sample location problem.  }

\section{Introduction}\label{sec:intro}

\subsection{David Tyler, beyond affine equivariance and elliptical symmetry}\label{sec:intro1}
The closely related concepts of affine equivariance and elliptical symmetry  played a central role in the development of robust multivariate statistics over the past 60 years.\footnote{\cite{Tukey60}, \cite{Huber64}, and \cite{Hampel68} generally are considered as   laying the foundations of modern robust statistics; see \cite{Ronch06} for a historical perspective and \cite{Stigler73} for an account of the pre-Tukey era.} A critical  attitude towards this dominant role of elliptical densities constitutes a red thread running through all of David's contributions to   multivariate anlysis\footnote{Significantly, {\it``Robust Multivariate Statistics: Beyond Ellipticity and Affine Equivariance''} is the title of one of David's NSF grants.}---an attitude that actually  takes  place in a broader debate on the ordering of the real space in dimension $d\geq 2$. Such  ordering  is an essential issue if the univariate concepts of distribution and quantile functions, ranks, and signs,   all heavily depending on the canonical ordering of the real line,  are to be extended to a multivariate context.  

\subsection{Ordering the real space in dimension $d\geq 2$}\label{sec:intro2}
The problem of ordering $\mathbb R ^d$ for $d\geq 2$, hence ranking multivariate observations,  has a long history in statistics. Many attempts have been made to define adequate  multivariate concepts of ranks. 

The notion of  ranks, however, is not an isolated one, as it is inseparable from that of empirical quantiles, quantile regions (collections of points with ranks less than or equal to  some given value), and quantile contours (collections of points with ranks  equal to some given value). A sound definition thus should include all these concepts, along with their population versions---the population distribution and quantile functions $F$ and $Q$---and their mutual relations (a  quantile function is the inverse of a distribution function $(Q=F^{-1})$; a population distribution function and its empirical version are asymptotically related via a Glivenko-Cantelli result, etc.).  Among the key properties of any successful concept are the distribution-freeness (within the class of absolutely continuous distributions $\rm P$, say) of the ranks and the distribution-freeness of the push-forward\footnote{We adopt here the convenient terminology and notation of measure transportation: the push-forward $F\#{\rm P}$ of $\rm P$ by $F$ is the distribution of $F(Z)$ where $Z\sim{\rm P}$, i.e.,  $F(Z)\sim F\#{\rm P}$ if $Z\sim\rm P$.}    of a distribution $\rm P$ by its distribution function $F$. Without these distribution-freeness properties,  the level of a quantile~$Q(\tau)=F^{-1}(\tau)$  depends on the distribution $\rm P$ characterized by $F$ and can be anything  larger or smaller than $\tau$: as a quantile of order $\tau$, thus, $Q(\tau)$ is totally meaningless.

Appealing as they are, none of the attempts that had been made until recently---marginal ranks, spatial ranks, elliptical (or Mahalanobis) ranks, ... ---is satisfying the desired properties; actually, none of them is  even enjoying distribution-freeness!\footnote{The Mahalanobis ranks and signs \citep{HP02} are enjoying distribution-freeness over the class of elliptically symmetric distributions only.} Nor do the various depth  concepts: the probability content of a depth contour of given depth strongly depends on the underlying $\rm P$, which hinders its interpretation as a quantile contour.

Based on measure transportation results (mainly, a theorem by \cite{McCann}) \cite{Chernoetal17}, \cite{Hallin17}, and \cite{Hallinetal21} recently introduced the  concepts  of {\it center-outward ranks} and {\it signs}, {\it distribution} and {\it quantile functions} which, for the first time, satisfy all the desired properties (see   \cite{Hallinetal21}  and the review by \cite{Hallin21} for details) and   further triggered the development of several  appealing multivariate, distribution-free statistical procedures, among which   \citet{FaugRusch17}, \citet{shi2020distribution,shi2020rate}, \citet{GhosSen19}, \citet{debsen21}, \citet{DebGhoSen20,debbhatsen21}, \citet{MasudAeron21}, \citet{HHH20,HLL19,HLL20}. 
 These are the concepts we are considering in this paper and now describe under their various versions.

\section{Center-outward ranks and signs}\label{sec:centeroutward}
For the simplicity of exposition, we throughout consider  distributions $\rm P$ on $\mathbb R ^d$ in the family $\mathcal P ^d$ of Lebesgue-absolutely continuous distributions with {\it nonvanishing densities},  that is,  with a density $f$ such that for all $B>0$ there exist $0<m^-_B\leq~\!m^+_B<~\!\infty$ such that $m^-_B\leq f(\zb)\leq m^+_B$ for all $\zb$ such that $\Vert\zb\Vert\leq B$. That assumption can be relaxed, though: see \cite{delb20}.

\subsection{A   measure transportation-based concept of distribution and quantile functions }\label{sec:centeroutward1}
The basic idea behind the definition of the center-outward distribution and quantile functions of a probability measure $\rm P\in\mathcal P ^d$ is quite simple. For $d=1$, the distribution function $F$ of $\rm P$ is the unique monotone increasing function pushing $\rm P$ forward to the uniform ${\rm U}_{[0,1]}$ over $[0,1]$---namely, $F\#{\rm P}={\rm U}_{[0,1]}$. Rather than $F$, however, which is based on a left-to-right ordering of $\mathbb R$ that does not extend to $\mathbb R ^d$ for $d\geq 2$, we consider the {\it center-outward distribution function} $F_\pm:=2F-1$, which contains the same information as $F$ and is the unique monotone increasing function pushing~$\rm P$ forward to the uniform ${\rm U}_{[-1,1]}$ over $[-1,1]$. A monotone increasing function is the gradient (the derivative) of a convex function: the center-outward distribution function $F_\pm$ this is the unique gradient of a convex function such that $F_\pm\#{\rm P}={\rm U}_{[-1,1]}$.  
The interval $[-1,1]$ 
 is, for $d=1$, the closed unit ball $\bar{\mathbb S}_d$, where ${\mathbb S}_d:=\big\{\ub \vert \Vert\ub\Vert \leq 1\big\}$ and, denoting by ${\rm U}_d$ the spherical uniform\footnote{The {\it spherical uniform} ${\rm U}_d$ over $\bar{\mathbb S}_d$ is the spherical distribution with center $\mathbf 0$ and  radial density the uniform over $[0,1]$: it is thus the product of a uniform over $[0,1]$ for the distances to the origin and a uniform over the unit (hyper)sphere for the directions.} over  $\bar{\mathbb S}_d$, the spherical uniform over $\bar{\mathbb S}_1$ coincides with the uniform over $[-1,1]$:  namely, ${\rm U}_1={\rm U}_{[-1,1]}$. 

A celebrated theorem by  \cite{McCann} tells us that, for arbitrary dimension~$d\in~\!\mathbb N$ and arbitrary~$\rm P\in\mathcal P ^d$, there exist a ($\rm P$-a.s., here Lebesgue-a.e.) unique gradient  of a convex function $\Fb_\pm$ such that $\Fb_\pm\#{\rm P}={\rm U}_d$. Obviously, for $d=1$, $\Fb_\pm$ coincides with the univariate $F_\pm$, whence the notation. Call $\Fb_\pm$ the {\it center-outward distribution function} of $\rm P$. It follows from \cite{Fig18} that---except perhaps at~$\Fb_\pm^{-1}({\mathbf 0})$ (a set of points with Lebesgue measure zero)---$\Fb_\pm\#{\rm P}$ is a homeomorphism, hence admits a continuous (except perhaps at $\Fb_\pm^{-1}({\mathbf 0})$) inverse $\Qb_\pm:=\Fb_\pm^{-1}$: call $\Qb_\pm$ the {\it center-outward quantile function} of $\rm P$. Clearly, $\Qb_\pm\#{\rm U}_d={\rm P}$.

This, with the spherical uniform ${\rm U}_d$ as a reference distribution (extending ${\rm U}_{[-1,1]}$)    is the concept proposed in  \cite{Hallin17} and \cite{Hallinetal21}, where  we refer to for further properties of $\Fb_\pm$ and $\Qb_\pm$ justifying their qualification as distribution and quantile functions.  

Other   choices are possible for the reference $\rm U$, though.  Replacing ${\rm U}_d$ with an arbitrary compactly supported  absolutely continuous reference distribution $\rm U$, \cite{Chernoetal17}, in a very general  approach, propose, under the name of {\it Monge-}{\it Kantorovich  vector rank} and  {\it Monge-}{\it Kantorovich  quantile} functions,  measure transportation-based definitions of a broad class of  analogues, $\Fb_{\text{\rm MK}}$ and~$\Qb_{\text{\rm MK}}$, say, of $\Fb_\pm$ and $\Qb_\pm$. For nonspherical $\rm U$'s, however, the Monge-Kantorovich  quantile functions do not enjoy all the features expected from a quantile function;\footnote{On this point, see Section 3.4 in \cite{Hallin21}.}  \cite{Chernoetal17} therefore also introduce a concept of {\it Monge-}{\it Kantorovich  depth} $D_{\text{\rm MK}}$---a transformation-retransformation version (based on the Monge-Kantorovich  vector rank function) of  classical Tukey depth $D_{\text{\rm Tukey}}$.  For spherical~$\rm U$'s, Monge-Kantorovich  depth and quantile contours coincide. More precisely, defining  $\delta (\tau):=D_{\text{\rm Tukey}}(\ub_\tau)$ where $\ub_\tau$ is such that ${\rm U}\big( \big\{ \ub\,\big\vert\, \Vert\ub\Vert\leq\Vert\ub_\tau\Vert\big\}\big) = \tau$, one has~$\Big\{\zb\, \big\vert\,  \Vert\Fb_{\text{\rm MK}}\Vert =~\!\tau\Big\}=~\!D_{\text{\rm MK}}^{-1}(\delta(\tau))$. Recurring to depth in order to construct quantile regions and contours, thus, is not necessary in the case of  a spherical reference~$\rm U$ which, in that respect, offers a better conceptual coherence between the resulting notions of vector ranks and quantiles. As far as rank tests are concerned, however, this can be considered a minor concern. 

The choice for $\rm U$ of the nonspherical  Lebesgue uniform ${\rm U}_{[0,1]^d}$ over the unit  (in the canonical basis) hypercube~$[0,1]^d$---call it the {\it cubic uniform}---yields a vector rank function $\Fb_{\text{\rm MK}}$  that  reduces, for $d=1$, to the classical   distribution   function~$F$ just as~$\Fb_\pm$ reduces to $F_\pm$. Despite poor equivariance properties\footnote{Contrary to $\Fb_\pm$, which is nicely equivariant, the rank vector function $\Fb_{\text{\rm MK}}$ associated with  the cubic uniform ${\rm U}_{[0,1]^d}$ is highly non-equivariant under orthogonal transformations.} its use has been advocated by several authors: see, e.g., \cite{FaugRusch17}, \cite{Carletal16}, 
\cite{DebGhoSen20,debbhatsen21}, \cite{debsen21}.
 \vspace{-1mm}

\subsection{Multivariate ranks and signs}\label{ranksignsec}

Denote by $\Zb\n:=(\Zb\n_1,\ldots,\Zb\n_n)$ an i.i.d.\ sample with distribution $\rm P \in \mathcal P ^d$. The empirical counterpart $\Fb_\pm\n$ of $\Fb_\pm$ is obtained as the solution of an optimal pairing problem between the sample values $\Zb\n_1,\ldots,\Zb\n_n$ and a ``regular'' grid ${\mathfrak G}\n$ with gridpoints ${\scriptstyle{\mathfrak{G}}}\n_1,\ldots,{\scriptstyle{\mathfrak{G}}}\n_n$. Precisely, $\left(\Fb_\pm\n(\Zb\n_1),\ldots, \Fb_\pm\n(\Zb\n_n)\right)$ is defined as the minimizer  $\left(
{\scriptstyle{\mathfrak{G}}}\n_{\pi^*(1)},\ldots,{\scriptstyle{\mathfrak{G}}}\n_{\pi^*(n)}\right)$, over the $n!$ possible permutations $\pi\in\Pi_n$ of the integers $\{1,\ldots,n\}$, of
$
  \sum_{i=1}^n\big\Vert {\bf Z}\n_{i} - {\scriptstyle{\mathfrak{G}}}\n_{\pi(i)})
\big\Vert ^2 
$.  

The choice of the grid ${\mathfrak G}\n$, of course, depends on the reference distribution $\rm U$ adopted in the definitions of Section~\ref{sec:centeroutward1}: in particular, the uniform discrete distribution over the $n$ gridpoints ${\scriptstyle{\mathfrak{G}}}\n_1,\ldots,{\scriptstyle{\mathfrak{G}}}\n_n$ should converge weakly to $\rm U$ as $n\to\infty$. Our objective is to investigate the finite-sample performance of two-sample location tests based on 
\begin{itemize}
\item[(Ti)] $\:\,\,$the empirical center-outward distribution function ${\bf F}\n_\pm$ associated with the spherical uniform reference distribution $\rm U =\rm U _d$;   
\item[(Tii)] $\:$the empirical Monge-Kantorovich vector ranks ${\bf F}\n_{\scriptscriptstyle{\square}}$ associated  with the cubic uniform  refe\-rence distribution~$\rm U =\rm U _{[0,1]^d}$;
\item[(Tiii)] $\:$the empirical Monge-Kantorovich vector ranks ${\bf F}\n_{\pm{\cal N}}$ associated  with the Gaussian~${\mathcal N}({\mathbf 0}, {\mathbf I}_d)$ reference considered as a spherical distribution; 
 \item[(Tiv)] $\:$the empirical Monge-Kantorovich vector ranks ${\bf F}\n_{\scriptscriptstyle{\square}{\cal N}}$ associated  with the Gaussian~${\mathcal N}({\mathbf 0}, {\mathbf I}_d)$ reference considered as  a product of univariate standard normal distributions.
\end{itemize}
The grids we are using for these four cases are constructed as follows (see Section~\ref{Haltonsec} for details on Halton sequences and the choice of $n_R$ and~$n_S$):
\begin{itemize}
\item[($\mathfrak G$i)] $\:\, \rm U = \rm U _d$: (a) factorize $n$ into $n=n_Rn_S + n_0$ with $n_0<\min(n_R,n_S)$; (b) generate a Halton sequence ${\mathfrak S}_{(n_S)}:=({\bf u}_1,\ldots {\bf u}_{n_S})$ over the unit (hyper)\-sphere~${\mathcal S _{d-1}}$; (c)~the grid ${\mathfrak G}\n$ consists of the intersections of these $n_S$ unit vectors with   the $n_R$ hyperspheres centered at $\mathbf 0$, with radii ${j}/{(n_R  +1)}$, 
$j=1,\ldots ,n_R$---along with $n_0$ copies of the origin;
\item[($\mathfrak G$ii)] $\:\!\! \rm U = \rm  U _{[0,1]^d}$: the grid ${\mathfrak G}\n$ is a Halton sequence over $[0,1]^d$;
\item[($\mathfrak G$iii)] $\:\!\!\!\hspace{-.2mm} \rm U ={\mathcal N}({\mathbf 0}, {\mathbf I}_d)$, spherical grid: the grid ${\mathfrak G}\n$ is the image, by  the radial transformation $\zb\mapsto \sqrt{F_{\chi^2_d}^{-1}(\Vert\zb\Vert)}\zb$,  of the spherical grid constructed in (i), where $F_{\chi^2_d}$ denotes the chi-square distribution function with $d$ degrees of freedom; 
\item[($\mathfrak G$iv)] $\:\!\!\! \rm U ={\mathcal N}({\mathbf 0}, {\mathbf I}_d)$, cubic grid: the grid ${\mathfrak G}\n$ is the image, by componentwise application of the standard normal quantile transformation $z_i\mapsto \Phi^{-1}(z_i)$,\footnote{As usual, we denote by $\Phi$ the standard normal distribution function, by $\Phi^{-1}$ the standard normal quantile function.} 
of a Halton sequence over  $[0,1]^d$.
\end{itemize}

\begin{remark} The grid ${\mathfrak G}\n$ in ($\mathfrak G$i) reduces, for $d=1$, to 
\[
\left\{\pm 1/(\lceil n/2\rceil +1), \ldots, \pm \lceil n/2\rceil/(\lceil n/2\rceil +1)\right\} 
\]
 along with the origin in case $n$ is odd; that grid is of the form 
\[
 \left\{2\left(1/(n+1)
\right)-1,\ldots, 2\left(n/(n+1)
\right)-1
\right\}
\]
where $ \left\{\left(1/(n+1)
\right),\ldots, \left(n/(n+1)
\right)
\right\}
$ is the grid producing traditional univariate ranks to which the grid ${\mathfrak G}\n$ in ($\mathfrak G$ii) also reduces for $d=1$. 
\end{remark}

\begin{remark} In ($\mathfrak G$i) and ($\mathfrak G$iii), the grid ${\mathfrak G}\n$ is  spherical; as a consequence, ${\bf F}\n_\pm(\Zb\n_i)$ and ${\bf F}\n_{\pm {\cal N}}(\Zb\n_i)$ in (Ti) and (Tiii) naturally factorize as\vspace{-1mm}   
$$ 
{\bf F}\n_\pm(\Zb\n_i)=\Vert{\bf F}\n_\pm(\Zb\n_i) \Vert\frac{{\bf F}\n_\pm(\Zb\n_i)}{\Vert{\bf F}\n_\pm(\Zb\n_i) \Vert}=: \frac{R\n_{i\pm}}{n_R+1}{\Sb\n_{i\pm}}
\vspace{-1mm} $$
where  $R\n_{i\pm}= (n_R+1)\Vert{\bf F}\n_\pm(\Zb\n_i)\Vert$, ranging from 0 or 1 (according as $n_0\neq 0$ or~$n_0=0$) to $n_R$, is the {\it center-outward rank} of $\Zb\n_i$ and $\Sb\n_{i\pm}$ (a unit vector) has the interpretation of a (multivariate) {\it center-outward sign} and\vspace{-1mm} 
\begin{equation}\label{vdWDef}
{\bf F}\n_{\pm {\cal N}}(\Zb\n_i)=\Vert{\bf F}\n_{\pm {\cal N}}(\Zb\n_i) \Vert
 \frac{{\bf F}\n_{\pm {\cal N}}(\Zb\n_i)}{
 \Vert{\bf F}\n_{\pm {\cal N}}(\Zb\n_i) \Vert
 }=: J_{\text{\rm vdW}}\left(\frac{R\n_{i\pm{\cal N}}}{n_R+1}\right){\Sb\n_{i\pm{\cal N}}}
\vspace{-1mm} \end{equation}
where $J_{\text{\rm vdW}}= \sqrt{F_{\chi^2_d}^{-1}}$ is the univariate {\it normal} or {\it van der Waerden score function},~${R\n_{i\pm{\cal N}}}$ the rank of~$\Vert{\bf F}\n_{\pm {\cal N}}(\Zb\n_i) \Vert$ among the $n_R$ distinct values of $\Vert{\bf F}\n_{\pm {\cal N}}(\Zb\n_i) \Vert$ for~$i=1,\ldots,n$ and ${\Sb\n_{i\pm{\cal N}}}$ similarly has the interpretation of a multivariate   sign. Being based on different transport maps, however, neither $R\n_{i\pm}$ and ${R\n_{i\pm{\cal N}}}$ nor $\Sb\n_{i\pm}$ and ${\Sb\n_{i\pm{\cal N}}}$ need coincide.
\end{remark}

\begin{remark} No similar factorization into ranks and signs occurs with the vector ranks~${\bf F}\n_{\scriptscriptstyle{\square}}$ and ${\bf F}\n_{\scriptscriptstyle{\square}{\cal N}}$ in (Tii) and (Tiv).\vspace{-2mm} 
\end{remark}

\subsection{Distribution-free tests based on center-outward ranks and signs }\label{sec:centeroutward2}
\cite{HHH20} propose, for multiple-output regression models with unspecified noise distribution $\rm P\in \mathcal P ^d$,  fully distribution-free yet, for adequate choice of   scores, parametrically efficient center-outward rank tests of the null hypothesis of no-treatment effect based on the empirical center-outward distribution functions $\Fb\n_\pm$ (hence, the center-outward ranks and signs). 

The particular case of  two-sample location is treated by \cite{debbhatsen21} who also consider tests based on the empirical Monge-Kantorovich vector ranks   ${\bf F}\n_{\text{\rm MK}}$ associated with various reference distributions. 

\subsubsection{Score functions
}

In line with the classical theory developed, e.g., by \cite{HajekSidak}, rank-based statistics, irrespective of the reference distribution, involve {\it score functions} or {\it scores.}  Depending on the context, a   score function is  a   mapping $\Jb$ from the unit ball ${\mathbb S}_d$ or the unit cube $[0,1]^d$ 
 to ${\mathbb R}^d$ satisfying some mild regularity assumptions    (continuity, square integrability, etc.: see, e.g. \cite{HHH20}, Assumption~3.1). The only score functions we are considering here are the {\it Wilcoxon}, the {\it spherical  van der Waerden}, and the {\it marginal van der Waerden  score functions}
 $$\Jb_{\text{\rm W}} (\ub)\!:=\!\ub,\quad\!\! \Jb^\pm_{\text{\rm vdW}}(\ub)\!:=\!\sqrt{F_{\chi^2_d}^{-1}}(\Vert\ub\Vert)\frac{\ub}{\Vert\ub\Vert},\text{ and } \Jb^{\scriptscriptstyle\square}_{\text{\rm vdW}}(\ub)\!:=\!\left(\Phi^{-1}(u_1),\ldots,\Phi^{-1}(u_d)
 \right),
 $$
 where $F_{\chi^2_d}$ and $\Phi$  stand for the (univariate) chi-square ($d$ degrees of freedom) and standard normal distribution functions, respectively.\vspace{-1mm}   

\subsubsection{Test statistics}\label{sec:teststat}

 For simplicity, our investigation here
  is limited to that particular case of   two-sample location models, where   the $n$ observations   are i.i.d.\ under the null and  consist of two samples,~$\Zb\n_1,\ldots,\Zb\n_{n_1}$ and $\Zb\n_{n_1+1},\ldots,\Zb\n_{n_1+n_2}$, with $n_1+n_2=n$. The classical procedure for this problem is Hotelling's test, based on  a quadratic statistic of the~form 
 $$\left(T _{\text{\rm Hot}}\n\right)^2= {\Deltab}_{\text{\rm Hot}}^{(n)\prime}\left({\boldsymbol\Sigma}\n_{\text{\rm Hot}}\right)^{-1} {\Deltab}\n_{\text{\rm Hot}}
 $$
 where ${\boldsymbol\Sigma}\n_{\text{\rm Hot}}$ is the estimated (under the null) covariance matrix of\vspace{-1mm}  
 \[
 {\Deltab}\n_{\text{\rm Hot}}:=\frac{1}{n_1}\sum_{i=1}^{n_1}\Zb\n_i - \frac{1}{n_2}\sum_{i=n_1+1}^{n}\Zb\n_i .
\vspace{-1mm}  \]
  The Hotelling test is parametrically efficient under Gaussian assumptions; it remains asymptotically valid,\footnote{Asymptotically valid, here, means pointwise (with respect to the actual density of the observations) asymptotically correct nominal probability  levels, not {\it uniformly}  asymptotically correct nominal probability  levels.} however, under mild moment assumptions and therefore qualifies as a pseudo-Gaussian procedure. 

 For score functions $\Jb$, the center-outward rank-based test statistics in Section~5.3.1 of~\cite{HHH20}  
are of the form\vspace{-1mm}  
\begin{equation}
\label{2samplestatpm}
\left(\tenq{T}^{(n)}_{\Jb\pm}\right)^2=\tenq{\Deltab}_{\Jb\pm}^{(n)\prime}\left({\boldsymbol\Sigma}_{\scriptscriptstyle\tenq{\scriptstyle\Deltab}_{\Jb\pm}}\right)^{-1}\tenq{\Deltab}\n_{\Jb\pm}
\end{equation}
where ${\boldsymbol\Sigma}_{\scriptscriptstyle\tenq{\scriptstyle\Deltab}_\Jb}$ is the exact or asymptotic covariance of \vspace{-1mm}
\begin{equation}
\label{2sampledelta}
\tenq{\Deltab}\n_{\Jb\pm}:=\frac{1}{n_1}\sum_{i=1}^{n_1} \Jb(\Fb\n_\pm(\Zb\n_i)) 
-\frac{1}{n_2}\sum_{i=n_1+1}^{n} \Jb(\Fb\n_\pm(\Zb\n_i)) .
\vspace{-1mm} 
\end{equation}
Since the quadratic form~\eqref{2samplestatpm} is invariant under affine  transformations of $\tenq{\Deltab}\n_{\Jb\pm}$ and since the sum~$\sum _{i=1}^n\Jb(\Fb\n_\pm(\Zb\n_i))$ is a deterministic constant that only depends on $\Jb$ and the grid used in the definition of $\Fb\n_\pm$, the same test statistic can be based on\vspace{-1mm}  
$$\tenq{\Deltab}^{(n)}_\Jb= \frac{1}{n_1} \sum_{i=1}^{n_1} \Jb(\Fb\n_\pm(\Zb\n_i)) -\frac{1}{n}\sum_{i=1}^n  \Jb(\Fb\n_\pm(\Zb\n_i)),
\vspace{-1mm} $$
yielding the test statistic described in Section~5.3.1 of~\cite{HHH20} which, in the particular case of  Wilcoxon and van der Waerden scores $\Jb_{\text{\rm W}} $ and  $\Jb^\pm_{\text{\rm vdW}} $,  we denote as $\left(\tenq{T}^{(n)}_{\text{\rm W}\pm}\right)^2$ and $\left(\tenq{T}^{(n)}_{\text{\rm vdW}\pm}\right)^2$, respectively.\smallskip

For the same testing problem,  \cite{debbhatsen21} consider statistics of the form \eqref{2samplestatpm}, but also \begin{compactenum}
\item[(a)] based on the empirical Monge-Kantorovich vector ranks ${\bf F}\n_{\scriptscriptstyle\square}
$ associated  with the cubic uniform  refe\-rence~$\rm U ={\rm U}_{[0,1]^d}$, statistics $\tenq{T}^{(n)}_{\Jb{\scriptscriptstyle\square}}$ and $\tenq{\Deltab}\n_{\Jb{\scriptscriptstyle\square}}$ of the same form as~$\tenq{T}^{(n)}_{\Jb\pm}$ and $\tenq{\Deltab}\n_{\Jb\pm}$ in~\eqref{2samplestatpm} and~\eqref{2sampledelta} but  with $\Jb(\Fb\n_{\scriptscriptstyle\square}(\Zb\n_i)) $ substituting $\Jb(\Fb\n_\pm(\Zb\n_i))$;  
denote by $\left(\tenq{T}^{(n)}_{\text{\rm W}{\scriptscriptstyle\square}}\right)^2$ and $\left(\tenq{T}^{(n)}_{\text{\rm vdW}{\scriptscriptstyle\square}}\right)^2$  the particular cases of  the Wilcoxon and cubic van der Waerden statistics obtained for the scores $\Jb_{\text{\rm W}} $ and  $\Jb^{\scriptscriptstyle\square}_{\text{\rm vdW}}$, respectively;
 
\item[(b)] based on the empirical Monge-Kantorovich vector ranks ${\bf F}\n_{\pm{\mathcal N}}
$ and ${\bf F}\n_{\scriptscriptstyle\square{\mathcal N}}$ associated with the spherical Gaussian reference ${\mathcal N}({\mathbf 0}, {\mathbf I}_d)$ considered as spherical or as a product of independent uniforms,  statistics $\tenq{T}^{(n)}_{\Jb{\pm}{\mathcal N}}$ and $\tenq{T}^{(n)}_{\Jb{\scriptscriptstyle\square}{\mathcal N}}$   of the same form as~$\tenq{T}^{(n)}_{\Jb\pm}$ and $\tenq{\Deltab}\n_{\Jb\pm}$ in~\eqref{2samplestatpm} and~\eqref{2sampledelta} but  with $\Jb(\Fb\n_{\pm{\mathcal N}}(\Zb\n_i)) $ and $\Jb(\Fb\n_{\scriptscriptstyle\square{\mathcal N}}(\Zb\n_i)) $,  respectively,  substituting $\Jb(\Fb\n_\pm(\Zb\n_i))$; this, for Wilcoxon scores $\Jb(\ub)=\ub\vspace{-1mm}$, yields the van der Waerden statistics $\left(\tenq{T}^{(n)}_{\text{\rm vdW}{\pm}{\mathcal N}}\right)^2$ and $\left(\tenq{T}^{(n)}_{\text{\rm vdW}{\scriptscriptstyle\square}{\mathcal N}}\right)^2$.
\end{compactenum}

\begin{remark} 
Although they are based on Wilcoxon (identity) scores, the terminology ``van der Waerden statistic'' for $\tenq{T}^{(n)}_{\text{\rm vdW}{\pm}{\mathcal N}}\vspace{1mm}$ and $\tenq{T}^{(n)}_{\text{\rm vdW}{\scriptscriptstyle\square}{\mathcal N}}$ seems more appropriate than the terminology ``Wilcoxon statistic''  used by \cite{debbhatsen21}, and is   in line with the traditional terminology of rank-based inference. Both $\tenq{T}^{(n)}_{\text{\rm vdW}{\pm}}\vspace{1mm}$ and $\tenq{T}^{(n)}_{\text{\rm vdW}{\pm}{\mathcal N}}$ indeed result from a transport from the sample values to a grid of Gaussian quantiles of the form ($\mathfrak G$iii). For   $\tenq{T}^{(n)}_{\text{\rm vdW}{\pm}}$, the transport is $\Jb\circ \Fb\n_\pm$, which, as a rule, is not an optimal one (not the gradient of a convex function) while, for  $\tenq{T}^{(n)}_{\text{\rm vdW}{\pm}{\mathcal N}}$, the transport is the optimal one ${\bf F}\n_{\scriptscriptstyle{\square}{\cal N}}$;  the difference between $\tenq{T}^{(n)}_{\text{\rm vdW}{\pm}}\vspace{1mm}$ and $\tenq{T}^{(n)}_{\text{\rm vdW}\pm{\mathcal N}}$    thus essentially consists in the way the transport to the spherical Gaussian grid is performed.  A similar remark can be made for $\tenq{T}^{(n)}_{\text{\rm vdW}{\scriptscriptstyle\square}}\vspace{1mm}$ and $\tenq{T}^{(n)}_{\text{\rm vdW}{\scriptscriptstyle\square}{\mathcal N}}$.\vspace{-2mm} 
\end{remark}



\section{Finite-sample performance: two-sample location simulations}
It clearly appears that  choices are to be made before performing a rank test based on the concepts  of multivariate ranks developed on the previous sections: center-outward ranks?  vector ranks? which ones? with which scores? The analysis of asymptotic performance does not help much, as the same local powers can be achieved by several  alternatives. The objective of this paper is to determine whether finite-sample performance can help us with these choices. We restrict ourselves to the two-sample location problem, Wilcoxon and van der Waerden scores, but the conclusions are quite likely to hold for the general case of multiple-output linear models considered in \cite{HHH20}. 


Before explaining how  simulations were conducted, let us provide some details on the way the grids described in Section~\ref{ranksignsec} were constructed. Recall that the aim of these grids is to provide a discrete approximation of the chosen continuous reference distribution.\vspace{-2mm} 

\subsection{Halton sequences on the cube and the sphere (($\mathfrak G$ii) and ($\mathfrak G$iv)  grids)\vspace{-1mm} }\label{Haltonsec}
The grid constructions ($\mathfrak G$ii) and ($\mathfrak G$iv) involve Halton sequences on the hypercube~$[0,1]^d$.  Halton sequences are pseudo-random numbers with low discrepancy of classical use in methods such as Monte Carlo simulations.  We used the implementation available in the package \textsf{SDraw}
by \citet{mcdonald2020package}.
The grid construction in ($\mathfrak G$i), hence also in ($\mathfrak G$iii), requires a $n_S$-point ``Halton sequence'' over the hypersphere~${\mathcal S}^{d-1}$. To obtain such a grid, we first generate a $n_S$-point Halton sequence over $[0,1]^{d-1}$ then componentwise perform the standard normal quantile transformation $u_j\mapsto z_j:=\Phi^{-1}(u_j)$. This yields $n_S$ points $\zb_1,\ldots,\zb_{n_S}$, with\vspace{-1mm} 
$$\zb_j:= (\Phi^{-1}(u_{j1}),\ldots, \Phi^{-1}(u_{jd})).$$
 The   resulting unit vectors $\zb_j/\Vert\zb_j \Vert$, $j=1,\ldots,n_S$ constitute the desired sequence over~${\mathcal S}^{d-1}$.\vspace{-2mm}\footnote{The justification is the fact that if the distribution of $\Zb$ is a product of independent univariate standard normal marginals, then $\Zb$ is spherical Gaussian ${\mathcal N}({\boldsymbol 0}, {\bf I})$ and hence $\Zb/\Vert\Zb\Vert$ is uniform over~${\mathcal S}^{d-1}$.%
  } 
 
\subsection{Factorization of $n$ (($\mathfrak G$i) and ($\mathfrak G$iii) grids)}\label{factsec} 
As for the grid constructions ($\mathfrak G$i) and ($\mathfrak G$iii), they require a factorization of $n$ into~$n_Rn_S +~\!n_0$ with~$n_0 <\min(n_R, n_S)$. Intuition suggests choosing $n_R$ and $n_S$ of  order $n^{1/d}$ and~$n^{(d-1)/d}$, respectively. This, however, is of little help for finite $n$. Since the grid is supposed to provide an approximation of the spherical uniform, we rather proceed by minimizing the Wasserstein distance to the spherical uniform as proposed in \cite{mordant21}. 
 %
%
%
%
More precisely, considering the grid with $n_R$ radial points described in ($\mathfrak G$i), denote by  ${\rm G}\n_{n_R}$ the   discrete measure placing a probability mass $1/n$ on each of the~$n$ gridpoints    except for the origin which receives probability mass  $n_0/n$.  We adopt here the strategy suggested in \cite{mordant21}  by selecting the grid with $n_R^*$ radial points, where\vspace{-1mm} 
\begin{equation}\label{minfact}
n_R^*:=\argmin_{ 1\le n_R\le n} \ W_2( {\rm G}\n_{n_R},\mathrm{U}_d)\vspace{-1mm} 
\end{equation}
($W_2$, as usual,  stands for the Wasserstein distance of order two). For $d\geq 3$, that distance~$W_2( {\rm G}\n_{n_R},\mathrm{U}_d)$ does not only depend on $n_R$ (hence on $n_S$) but also on the~$n_S$ points chosen (as explained in Section~\ref{Haltonsec})   on the hypersphere ${\mathcal S}^{d-1}$.   The minimization in~\eqref{minfact} is feasible, as $n$, $n_R$, $n_S$, and $n_0$ all are integers. A similar strategy is adopted for the construction of the spherical Gaussian grids ($\mathfrak G$iii).  

Table~\ref{table:1} provides, for dimensions $d=2$ and $d=5$, various sample sizes,  
and  reference distributions  the spherical uniform (($\mathfrak G$i) grids) and the spherical Gaussian (($\mathfrak G$iii) grids),  the ``optimal values"  obtained via \eqref{minfact} for $n_R$, $n_S$, and $n_0$. These values are in line with the intuition that the ``optimal'' $n_R$ behaves like $n^{1/d}$ while the role of distances to the center rapidly decreases as the dimension $d$ increases.\vspace{-1mm}  
 
 \begin{table}[h!]
\centering
\begin{tabular}{|c | c | c c c c c |} 
 \hline \hline
 Reference distribution   & $d$  \ & $\quad n=50\quad $\ & $\quad n=100\quad $ \ & $\quad n=200\quad $ \ & $\quad n=300\quad $ \ & $\quad n=400\quad $\  \\ [0.5ex] 
 \hline\hline
 \vspace{-3mm}&  &  & & & &  \\ 
 &  & $n_R=4$ & $n_R=6$ & $n_R=9$ &$n_R=11$ &$n_R=12$ \\ 
 $\mathrm{U}_d$, ($\mathfrak G$i) grid &2 & $n_S=12$&$n_S=16$ & $n_S=22$ &$n_S=27$ &$n_S=33$ \\ 
  & &$n_0=2$&$n_0=4$ & $n_0=2$ &$n_0=3$ &$n_0=4$ \\ 
 \hline
 &  & $n_R=4$& $n_R=7$ & $n_R=11$ & $n_R=14$ & $n_R=18$ \\
  ${\mathcal N}({\mathbf 0}, {\mathbf I}_d)$,  ($\mathfrak G$iii) grid  &2 & $n_S=12$& $n_S=14$ & $n_S=18$ &$n_S=21$ &$n_S=22$ \\ 
 & & $n_0=2$& $n_0=2$ & $n_0=2$ &$n_0=6$ &$n_0=4$ \\  \hline 
  \vspace{-2mm}&  &  & & & &  \\ 
$n^{1/d}$ &  & $n^{1/2}\!=7.071$& $n^{1/2}\!=10.000$\vspace{1mm} & $n^{1/2}\!=14.142$ &$n^{1/2}\!=17.321$ &$n^{1/2}\!=20.000$ \\  %
 \hline
 \hline
\vspace{-3mm}&  & &  & & &  \\ 
  &  &  $n_R=2$&  $n_R=2$ &  $n_R=2$ &  $n_R=3$ &  $n_R=3$ \\
  $\mathrm{U}_d$, ($\mathfrak G$i) grid  & 5 &  $n_S=25$&  $n_S=50$ &  $n_S=100$ &  $n_S=100$ &  $n_S=133$ \\
  &  &  $n_0=0$&  $n_0=0$ &  $n_0=0$ &  $n_0=0$ &  $n_0=1$ \\
  \hline
 &  &   $n_R=1$&   $n_R=1$ &   $n_R=1$ &   $n_R=2$ &   $n_R=2$   \\
${\mathcal N}({\mathbf 0}, {\mathbf I}_d)$, ($\mathfrak G$iii) grid& 5 &   $n_S=50$&   $n_S=100$ &   $n_S=200$ &   $n_S=150$ &   $n_S=200$\\
  &  &   $n_0=0$&   $n_0=0$ &   $n_0=0$ &   $n_0=0$ &   $n_0=0$ \\   \hline
  \vspace{-2mm}&  &  & & & &  \\  
$n^{1/d}$ &  & $n^{1/5}\!=2.187$& $n^{1/5}\!=2.512$\vspace{1mm} & $n^{1/5}\!=2.885$ &$n^{1/5}\!=3.129$ &$n^{1/5}\!=3.314$ \\   \hline \hline
\end{tabular}
\caption{Optimal (in the sense of \eqref{minfact}) values of  $n_R$, $n_S$, and $n_0$ as  functions of the sample size $n$, the dimension $d$, and the reference distributions (($\mathfrak G$i) or ($\mathfrak G$iii) grids). \vspace{-5mm}}
\label{table:1}
\end{table}\vspace{-4mm}

\subsection{Simulations}
Based on the  grids obtained along the lines described in Sections~\ref{Haltonsec} and \ref{factsec},   the distribution-free critical values of the various rank tests under study were computed from 40 000 replications. Throughout, we chose $n_1=n_2=n/2$. The optimal maps between the sample and the grids were  obtained via an exact solver relying on the so-called Hungarian method that is implemented in the \textsf{R}-package \textsf{clue} by \citet{hornik2021package}. We now  turn to the empirical evaluation of the performance of the various rank-based Wilcoxon and van der Waerden  tests for the two-sample location problem. 

The objective of our simulations is, essentially,   empirical answers to the following two questions:
\begin{compactenum}
\item[(a)] should we use spherical grids (($\mathfrak G$i) or ($\mathfrak G$iii)) or cubic (($\mathfrak G$ii) or ($\mathfrak G$iv)) ones?
\item[(b)] should we, in line with the H\' ajek tradition,  privilege transports to the uniform combined  with scores  (as in $\left(\tenq{T}^{(n)}_{\text{\rm vdW}{\pm}}\right)^2$ and $\left(\tenq{T}^{(n)}_{\text{\rm vdW}{\scriptscriptstyle\square}}\right)^2$), or, as recommended by \cite{debbhatsen21}, should we rather consider ``direct transports''  to the ``scored distribution,'' that is, choose as reference distribution  the push-forward of the uniform by the score as in $\left(\tenq{T}^{(n)}_{\text{\rm vdW}{\pm}{\mathcal N}}\right)^2$ and $\left(\tenq{T}^{(n)}_{\text{\rm vdW}{\scriptscriptstyle\square}{\mathcal N}}\right)^2$?    
\end{compactenum}
\vspace{-2mm} 

%
%

\section{Wilcoxon-type  tests}
\label{sec: Wilcox}

The Wilcoxon tests are based on the identity score function $\Jb (\ub)=\ub$ and uniform (either spherical or cubic) reference distributions, yielding (see Section~\ref{sec:teststat}) the test statistics $\left(\tenq{T}^{(n)}_{\text{\rm W}{\pm}}\right)^2$ and~$\left(\tenq{T}^{(n)}_{\text{\rm W}{\scriptscriptstyle\square}}\right)^2$.\vspace{-3mm}  

%

\subsection{The bivariate case}
\label{sec: Wilcox2D}
In this section, we evaluate the performance of the  Wilcoxon tests based on~$\left(\tenq{T}^{(n)}_{\text{\rm W}{\pm}}\right)^2$ and~$\left(\tenq{T}^{(n)}_{\text{\rm W}{\scriptscriptstyle\square}}\right)^2\vspace{1mm}$ for samples of size~$n_1=n_2=n/2$ with $n=100$, 200, and~400. The first sample is drawn from a centered distribution, the second one  from a shifted version with shift $(\eta, \eta)\pr$, $\eta>0$ of the same.  The number of replications is $N=500$.
\vspace{-4mm}

\subsubsection{Spherical Gaussian samples}\label{sec411}
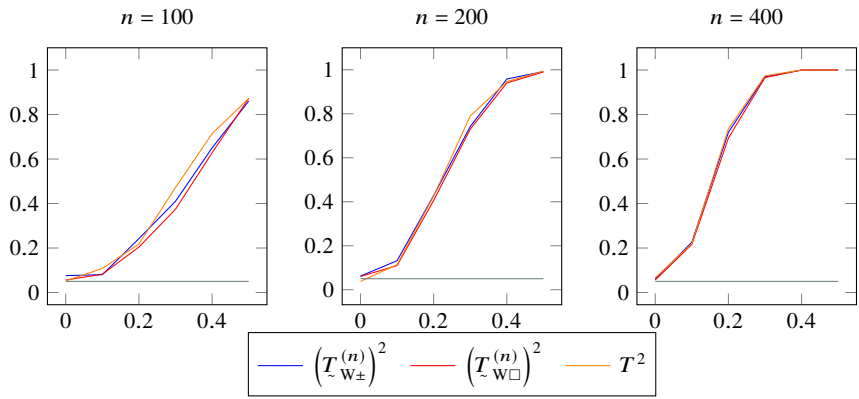
\begin{figure}[h]
\centering
\begin{tikzpicture}
    \begin{groupplot}[group style={
                      group name=myplot,
                      group size= 3 by 1
                      },height=5cm,width=4.5cm,
                      ytick={0,0.2,0.4,0.6,0.8,1},
                      ymax=1.1]
\nextgroupplot[title={$n=100$}]\addplot [domain=0:.5, color=blue] coordinates { 
(0,0.0760000000000001)
(0.1,0.0800000000000001)
(0.2,0.242)
(0.3,0.41)
(0.4,0.65)
(0.5,0.860000000000001)
}; 
\label{plots:plot1} 
\addplot [domain=0:.5, color=red] coordinates { 
(0,0.056)
(0.1,0.0820000000000001)
(0.2,0.204)
(0.3,0.374)
(0.4,0.628)
(0.5,0.872000000000001)
}; 
\label{plots:plot2} 
\addplot [domain=0:.5, color=orange] coordinates { 
(0, 0.052)
(0.1, 0.108)
(0.2, 0.218)
(0.3, 0.472)
(0.4, 0.714)
(0.5, 0.872)
}; 
\label{plots:plotH} 
\addplot [domain=0:.5, color=grey] coordinates { 
(0,0.05)
(0.1,0.05)
(0.2,0.05)
(0.3,0.05)
(0.4,0.05)
(0.5,0.05)
}; 
\label{plots:plot3} 
\nextgroupplot[title={$n=200$}]\addplot [domain=0:.5, color=blue] coordinates { 
(0,0.062)
(0.1,0.132)
(0.2,0.428)
(0.3,0.742000000000001)
(0.4,0.958000000000001)
(0.5,0.992000000000001)
}; 
\addplot [domain=0:.5, color=red] coordinates { 
(0,0.06)
(0.1,0.11)
(0.2,0.406)
(0.3,0.730000000000001)
(0.4,0.940000000000001)
(0.5,0.990000000000001)
}; 
\addplot [domain=0:.5, color=orange] coordinates { 
(0, 0.038)
(0.1, 0.114)
(0.2, 0.426)
(0.3, 0.790)
(0.4, .946)
(0.5, 0.994)
}; 

\addplot [domain=0:.5, color=grey] coordinates { 
(0,0.05)
(0.1,0.05)
(0.2,0.05)
(0.3,0.05)
(0.4,0.05)
(0.5,0.05)
}; 
\nextgroupplot[title={$n=400$}]\addplot [domain=0:.5, color=blue] coordinates { 
(0,0.06)
(0.1,0.226)
(0.2,0.720000000000001)
(0.3,0.970000000000001)
(0.4,1)
(0.5,1)
}; 
\addplot [domain=0:.5, color=red] coordinates { 
(0,0.056)
(0.1,0.216)
(0.2,0.694000000000001)
(0.3,0.966000000000001)
(0.4,1)
(0.5,1)
}; 
\addplot [domain=0:.5, color=orange] coordinates { 
(0, 0.066)
(0.1, 0.220)
(0.2, 0.736)
(0.3, 0.974)
(0.4, 1)
(0.5, 1)
}; 
\addplot [domain=0:.5, color=grey] coordinates { 
(0,0.05)
(0.1,0.05)
(0.2,0.05)
(0.3,0.05)
(0.4,0.05)
(0.5,0.05)
}; 

            \end{groupplot}      
\path (myplot c1r1.south west|-current bounding box.south)--
      coordinate(legendpos)
      (myplot c3r1.south east|-current bounding box.south);
\matrix[
    matrix of nodes,
    anchor=south,
    draw,
    inner sep=0.2em,
    draw
  ]at([yshift=-6ex]legendpos)
  {
    \ref{plots:plot1}& $\left(\tenq{T}^{(n)}_{\text{\rm W}{\pm}}\right)^2$ &[5pt]
    \ref{plots:plot2}& $\left(\tenq{T}^{(n)}_{\text{\rm W}{\scriptscriptstyle\square}}\right)^2$ &[5pt]
    \ref{plots:plotH}  & $T^2$ \\};
\end{tikzpicture}
\caption{\small\slshape Rejection frequencies, for spherical Gaussian samples (see 4.4.1) and various sample sizes, of Hotelling's test based on $T^2$ and the Wilcoxon tests based on $\left(\tenq{T}^{(n)}_{\text{\rm W}{\pm}}\right)^2$ and $\left(\tenq{T}^{(n)}_{\text{\rm W}{\scriptscriptstyle\square}}\right)^2$, respectively,  as  functions of the shift $\eta$;  
$N=500$ replications.\vspace{-8mm}} 
\label{fig: BivCompGausWilcox}
\end{figure}

\vspace{-3mm}
The first sample is drawn from~$\mathcal{N} ( (0, 0)\pr ,  {\bf I}_2)$, the second one from~$\mathcal{N} ( (\eta, \eta)^\top ,  {\bf I}_2)$ with~$\eta>0$.   Rejection frequencies over $N=500$ replications are
 shown (as  functions of $\eta$) in Figure~\ref{fig: BivCompGausWilcox}. All three tests display, essentially, the same performance: although Wilcoxon, in principle, is strictly less powerful than Hotelling (which in this case is finite-sample optimal), no significant loss of efficiency is detected.  
\vspace{-6mm}

\subsubsection{Spherical Student samples}\label{sec412}\vspace{-3mm}

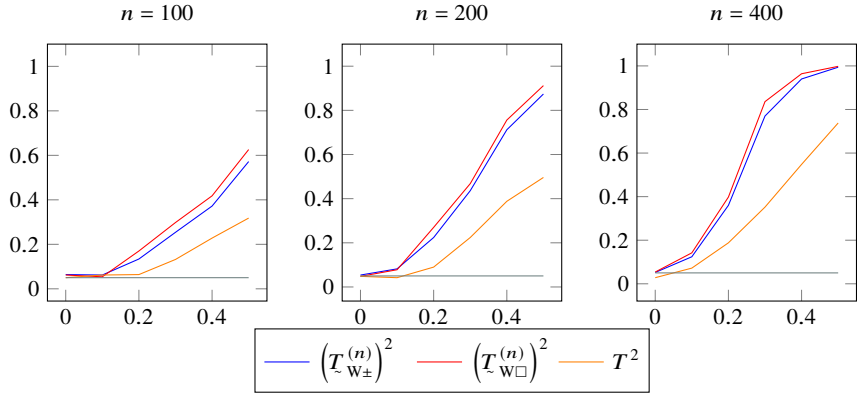
\begin{figure}[h]
\centering
\begin{tikzpicture}
    \begin{groupplot}[group style={
                      group name=myplot,
                      group size= 3 by 1
                      },height=5cm,width=4.5cm,
                      ytick={0,0.2,0.4,0.6,0.8,1},
                      ymax=1.1]
\nextgroupplot[title={$n=100$}]\addplot [domain=0:.5, color=blue] coordinates { 
(0,0.064)
(0.1,0.062)
(0.2,0.134)
(0.3,0.254)
(0.4,0.372)
(0.5,0.572)
}; 
\label{plots:plot1} 
\addplot [domain=0:.5, color=red] coordinates { 
(0,0.062)
(0.1,0.054)
(0.2,0.17)
(0.3,0.298)
(0.4,0.418)
(0.5,0.626)
}; 
\label{plots:plot2} 
\addplot [domain=0:.5, color=orange] coordinates { 
(0, 0.05)
(0.1, 0.062)
(0.2, 0.064)
(0.3, 0.132)
(0.4, 0.228)
(0.5, 0.318)
}; 
\label{plots:plotH} 
\addplot [domain=0:.5, color=grey] coordinates { 
(0,0.05)
(0.1,0.05)
(0.2,0.05)
(0.3,0.05)
(0.4,0.05)
(0.5,0.05)
}; 
\label{plots:plot3} 
\nextgroupplot[title={$n=200$}]\addplot [domain=0:.5, color=blue] coordinates { 
(0,0.054)
(0.1,0.0820000000000001)
(0.2,0.224)
(0.3,0.436)
(0.4,0.712000000000001)
(0.5,0.874000000000001)
}; 
\addplot [domain=0:.5, color=red] coordinates { 
(0,0.048)
(0.1,0.0780000000000001)
(0.2,0.27)
(0.3,0.468)
(0.4,0.756000000000001)
(0.5,0.912000000000001)
}; 
\addplot [domain=0:.5, color=orange] coordinates { 
(0, 0.048)
(0.1, 0.042)
(0.2, 0.09)
(0.3, 0.224)
(0.4, 0.388)
(0.5, 0.496)
}; 
\addplot [domain=0:.5, color=grey] coordinates { 
(0,0.05)
(0.1,0.05)
(0.2,0.05)
(0.3,0.05)
(0.4,0.05)
(0.5,0.05)
}; 
\nextgroupplot[title={$n=400$}]\addplot [domain=0:.5, color=blue] coordinates { 
(0,0.052)
(0.1,0.124)
(0.2,0.36)
(0.3,0.770000000000001)
(0.4,0.940000000000001)
(0.5,0.994000000000001)
}; 
\addplot [domain=0:.5, color=red] coordinates { 
(0,0.054)
(0.1,0.142)
(0.2,0.398)
(0.3,0.836000000000001)
(0.4,0.964000000000001)
(0.5,0.998000000000001)
}; 
\addplot [domain=0:.5, color=orange] coordinates { 
(0, 0.028)
(0.1, 0.072)
(0.2, 0.188)
(0.3, 0.352)
(0.4, 0.548)
(0.5, 0.738)
}; 
\addplot [domain=0:.5, color=grey] coordinates { 
(0,0.05)
(0.1,0.05)
(0.2,0.05)
(0.3,0.05)
(0.4,0.05)
(0.5,0.05)
}; 

            \end{groupplot}      
\path (myplot c1r1.south west|-current bounding box.south)--
      coordinate(legendpos)
      (myplot c3r1.south east|-current bounding box.south);
\matrix[
    matrix of nodes,
    anchor=south,
    draw,
    inner sep=0.2em,
    draw
  ]at([yshift=-6ex]legendpos)
  {
      \ref{plots:plot1}& $\left(\tenq{T}^{(n)}_{\text{\rm W}{\pm}}\right)^2$ &[5pt]
    \ref{plots:plot2}& $\left(\tenq{T}^{(n)}_{\text{\rm W}{\scriptscriptstyle\square}}\right)^2$ &
    \ref{plots:plotH} & $T^2$
    \\};
\end{tikzpicture}
\caption{\small\slshape Rejection frequencies, for Student samples (see 4.1.2) and various sample sizes, of Hotelling's test based on $T^2$ and the Wilcoxon tests based on $\left(\tenq{T}^{(n)}_{\text{\rm W}{\pm}}\right)^2$ and $\left(\tenq{T}^{(n)}_{\text{\rm W}{\scriptscriptstyle\square}}\right)^2$, respectively,  as  functions of the shift $\eta$;  
$N=500$ replications.\vspace{-5mm}} 
\label{fig: BivCompTsphWilcox}
\end{figure}

The first sample is drawn from a centered spherical Student with 2.1 degrees of freedom $ t_{2.1} ( (0,0)\pr, {\bf I}_2)$,   the second one from the shifted version~$ t_{2.1} ((\eta, \eta)\pr , {\bf I}_2)$ of the same distribution.   Rejection frequencies over $N=500$ replications are
 shown (as  functions of $\eta$) in Figure~\ref{fig: BivCompTsphWilcox}. The Wilcoxon tests substantially outperform Hotelling, an advantage that does not disappear with increasing $n$. Although the underlying distribution is spherical,   very slight superiority of~$\left(\tenq{T}^{(n)}_{\text{\rm W}{\scriptscriptstyle\square}}\right)^2\vspace{1mm}$ over $\left(\tenq{T}^{(n)}_{\text{\rm W}{\pm}}\right)^2$.\vspace{-4mm}

\subsubsection{``Banana-shaped''  samples}\label{sec413}\vspace{-3mm}

\begin{figure}[h]
\centering
\begin{tikzpicture}
    \begin{groupplot}[group style={
                      group name=myplot,
                      group size= 3 by 1
                      },height=5cm,width=4.5cm,
                      ytick={0,0.2,0.4,0.6,0.8,1},
                      ymax=1.1]
\nextgroupplot[title={$n=100$}]\addplot [domain=0:.5, color=blue] coordinates { 
(0,0.052)
(0.1,0.198)
(0.2,0.622)
(0.3,0.906000000000001)
(0.4,0.992000000000001)
(0.5,0.998000000000001)
}; 
\label{plots:plot1} 
\addplot [domain=0:.5, color=red] coordinates { 
(0,0.054)
(0.1,0.162)
(0.2,0.532)
(0.3,0.846000000000001)
(0.4,0.972000000000001)
(0.5,0.994000000000001)
}; 
\label{plots:plot2} 
\addplot [domain=0:.5, color=orange] coordinates { 
(0, 0.072)
(0.1, .1)
(0.2, .25)
(0.3, .508)
(0.4, .758)
(0.5, .886)
}; 
\label{plots:plotH} 
\addplot [domain=0:.5, color=grey] coordinates { 
(0,0.05)
(0.1,0.05)
(0.2,0.05)
(0.3,0.05)
(0.4,0.05)
(0.5,0.05)
}; 
\label{plots:plot3} 
\nextgroupplot[title={$n=200$}]\addplot [domain=0:.5, color=blue] coordinates { 
(0,0.044)
(0.1,0.37)
(0.2,0.900000000000001)
(0.3,1)
(0.4,1)
(0.5,1)
}; 
\addplot [domain=0:.5, color=red] coordinates { 
(0,0.044)
(0.1,0.292)
(0.2,0.840000000000001)
(0.3,0.984000000000001)
(0.4,1)
(0.5,1)
}; 
\addplot [domain=0:.5, color=orange] coordinates { 
(0, 0.04)
(0.1, .156)
(0.2, .444)
(0.3, .816)
(0.4, .976)
(0.5, .998)
}; 
\addplot [domain=0:.5, color=grey] coordinates { 
(0,0.05)
(0.1,0.05)
(0.2,0.05)
(0.3,0.05)
(0.4,0.05)
(0.5,0.05)
}; 
\nextgroupplot[title={$n=400$}]\addplot [domain=0:.5, color=blue] coordinates { 
(0,0.05)
(0.1,0.67)
(0.2,0.996000000000001)
(0.3,1)
(0.4,1)
(0.5,1)
}; 
\addplot [domain=0:.5, color=red] coordinates { 
(0,0.058)
(0.1,0.55)
(0.2,0.990000000000001)
(0.3,1)
(0.4,1)
(0.5,1)
}; 
\addplot [domain=0:.5, color=orange] coordinates { 
(0, 0.046)
(0.1, .248)
(0.2, .768)
(0.3, .984)
(0.4, 1)
(0.5, 1)
}; 
\addplot [domain=0:.5, color=grey] coordinates { 
(0,0.05)
(0.1,0.05)
(0.2,0.05)
(0.3,0.05)
(0.4,0.05)
(0.5,0.05)
}; 

            \end{groupplot}      
\path (myplot c1r1.south west|-current bounding box.south)--
      coordinate(legendpos)
      (myplot c3r1.south east|-current bounding box.south);
\matrix[
    matrix of nodes,
    anchor=south,
    draw,
    inner sep=0.2em,
    draw
  ]at([yshift=-6ex]legendpos)
  {
      \ref{plots:plot1}& $\left(\tenq{T}^{(n)}_{\text{\rm W}{\pm}}\right)^2$ &[5pt]
    \ref{plots:plot2}& $\left(\tenq{T}^{(n)}_{\text{\rm W}{\scriptscriptstyle\square}}\right)^2$ &[5pt]
    \ref{plots:plotH}& $T^2$
    \\};\end{tikzpicture}
\caption{\small\slshape Rejection frequencies, for ``banana-shaped'' samples (see 4.1.3) and various sample sizes, of Hotelling's test based on $T^2$ and the Wilcoxon tests based on $\left(\tenq{T}^{(n)}_{\text{\rm W}{\pm}}\right)^2$ and $\left(\tenq{T}^{(n)}_{\text{\rm W}{\scriptscriptstyle\square}}\right)^2$, respectively,  as  functions of the shift $\eta$;  
$N=500$ replications.\vspace{-6mm}} 
\label{fig: BivCompBananaWilcox}
\end{figure}
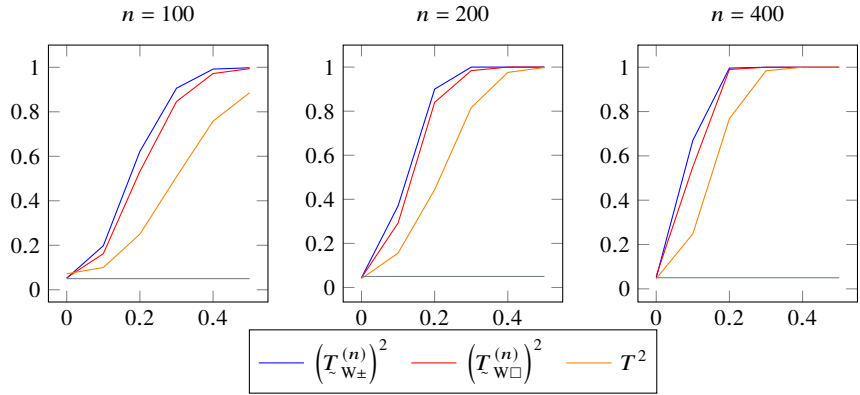

The first sample is drawn from a centered ``banana-shaped''   mixture 
\begin{align*}
\nonumber
0.3 
 \, \normal_2\left(
\begin{pmatrix} \phantom{-}0\phantom{.7} \\ -0.7 \end{pmatrix}, 
 \begin{pmatrix} 
 0.35^2&0 \\
0& 0.35^2 
 \end{pmatrix} 
  \right)&+ 
 0.35 
 \, \normal_2\left(
\begin{pmatrix} -0.9 \\ \phantom{-}0.3 \end{pmatrix}, 
 \begin{pmatrix} 
 \phantom{-}0.358&-0.55 \\
-0.55\phantom{0}& \phantom{-}1.02 
 \end{pmatrix} 
  \right)\\
\label{eq:banana}
  &+ 
  0.35 
   \, \normal_2\left(
\begin{pmatrix} \phantom{-}0.9 \\ \phantom{-}0.3 \end{pmatrix}, 
 \begin{pmatrix} 
 \phantom{-}0.358& \phantom{-}0.55 \\
\phantom{-}0.55\phantom{0}& \phantom{-}1.02 
 \end{pmatrix}
  \right).
 \end{align*}\vspace{-0mm}
of three Gaussian components. 
The second sample is drawn from a shifted version (shift $(\eta, \eta)\pr$, $\eta>0$) of the same mixture.   Rejection frequencies over $N=500$ replications are
 shown (as  functions of $\eta$) in Figure~\ref{fig: BivCompBananaWilcox}. The conclusions are the same  as in the previous case, except that the (slight) advantage now belongs to~$\left(\tenq{T}^{(n)}_{\text{\rm W}{\pm}}\right)^2$, despite the fact that the actual distribution is highly nonspherical.\vspace{-5mm}


\subsubsection{Samples   with independent Cauchy marginals}\label{sec414}\vspace{-3mm}

\begin{figure}[h]
\centering
\begin{tikzpicture}
    \begin{groupplot}[group style={
                      group name=myplot,
                      group size= 3 by 1
                      },height=5cm,width=4.5cm,
                      ytick={0,0.2,0.4,0.6,0.8,1},
                      ymax=1.1]
\nextgroupplot[title={$n=100$}]\addplot [domain=0:.5, color=blue] coordinates { 
(0,0.056)
(0.1,0.052)
(0.2,0.0880000000000001)
(0.3,0.116)
(0.4,0.194)
(0.5,0.29)
}; 
\label{plots:plot1} 
\addplot [domain=0:.5, color=red] coordinates { 
(0,0.046)
(0.1,0.062)
(0.2,0.0980000000000001)
(0.3,0.154)
(0.4,0.262)
(0.5,0.392)
}; 
\label{plots:plot2} 
\addplot [domain=0:.5, color=orange] coordinates { 
(0, 0.014)
(0.1, .012)
(0.2, .016)
(0.3, .026)
(0.4, .034)
(0.5, .028)
}; 
\label{plots:plotH} 
\addplot [domain=0:.5, color=grey] coordinates { 
(0,0.05)
(0.1,0.05)
(0.2,0.05)
(0.3,0.05)
(0.4,0.05)
(0.5,0.05)
}; 
\label{plots:plot3} 
\nextgroupplot[title={$n=200$}]\addplot [domain=0:.5, color=blue] coordinates { 
(0,0.054)
(0.1,0.0780000000000001)
(0.2,0.12)
(0.3,0.178)
(0.4,0.362)
(0.5,0.522)
}; 
\addplot [domain=0:.5, color=red] coordinates { 
(0,0.058)
(0.1,0.0800000000000001)
(0.2,0.15)
(0.3,0.24)
(0.4,0.464)
(0.5,0.652)
}; 
\addplot [domain=0:.5, color=orange] coordinates { 
(0, 0.018)
(0.1, .024)
(0.2, .022)
(0.3, .036)
(0.4, .046)
(0.5, .03)
}; 
\addplot [domain=0:.5, color=grey] coordinates { 
(0,0.05)
(0.1,0.05)
(0.2,0.05)
(0.3,0.05)
(0.4,0.05)
(0.5,0.05)
}; 
\nextgroupplot[title={$n=400$}]\addplot [domain=0:.5, color=blue] coordinates { 
(0,0.03)
(0.1,0.0760000000000001)
(0.2,0.206)
(0.3,0.398)
(0.4,0.622)
(0.5,0.786000000000001)
}; 
\addplot [domain=0:.5, color=red] coordinates { 
(0,0.03)
(0.1,0.0920000000000001)
(0.2,0.266)
(0.3,0.546)
(0.4,0.762000000000001)
(0.5,0.922000000000001)
}; 
\addplot [domain=0:.5, color=orange] coordinates { 
(0, 0.028)
(0.1, .026)
(0.2, .012)
(0.3, .014)
(0.4, .02)
(0.5, .038)
}; 
\addplot [domain=0:.5, color=grey] coordinates { 
(0,0.05)
(0.1,0.05)
(0.2,0.05)
(0.3,0.05)
(0.4,0.05)
(0.5,0.05)
}; 
            \end{groupplot}      
\path (myplot c1r1.south west|-current bounding box.south)--
      coordinate(legendpos)
      (myplot c3r1.south east|-current bounding box.south);
\matrix[
    matrix of nodes,
    anchor=south,
    draw,
    inner sep=0.2em,
    draw
  ]at([yshift=-6ex]legendpos)
  {
    \ref{plots:plot1}& $\left(\tenq{T}^{(n)}_{\text{\rm W}{\pm}}\right)^2$ &[5pt]
    \ref{plots:plot2}& $\left(\tenq{T}^{(n)}_{\text{\rm W}{\scriptscriptstyle\square}}\right)^2$ & [5pt]
    \ref{plots:plotH}& $T^2$
    \\};
    \end{tikzpicture}
\caption{\small\slshape Rejection frequencies, for samples with  independent Cauchy  marginals (see 4.1.4) and various sample sizes, of Hotelling's test based on $T^2$ and the Wilcoxon tests based on $\left(\tenq{T}^{(n)}_{\text{\rm W}{\pm}}\right)^2$ and~$\left(\tenq{T}^{(n)}_{\text{\rm W}{\scriptscriptstyle\square}}\right)^2$, respectively,  as  functions of the shift $\eta$;  
$N=500$ replications.\vspace{-7mm}} 
\label{fig: BivCompCindWilcox}
\end{figure}
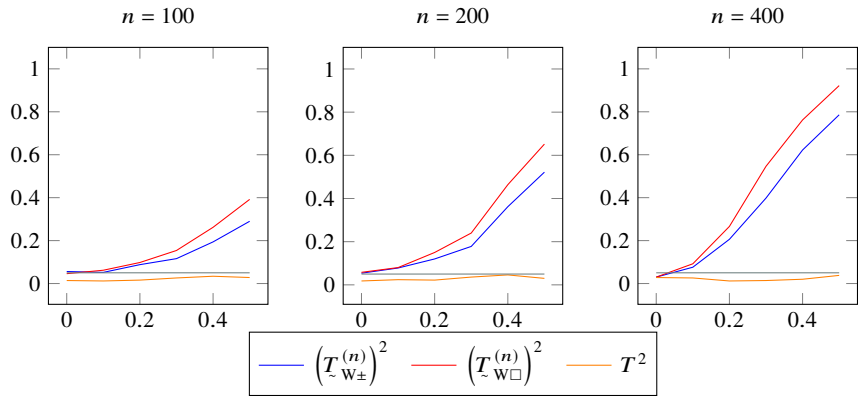

The first sample is drawn from a product of two independent Cauchy,   the second one from the shifted version of the same distribution.   Rejection frequencies over~$N=500$ replications are
 shown (as  functions of $\eta$) in  Figure~\ref{fig: BivCompCindWilcox}. With a rejection probability uniformly less than the nominal $5\%$ level, Hotelling, as expected, performs  miserably. In this independent component situation,  $\left(\tenq{T}^{(n)}_{\text{\rm W}{\scriptscriptstyle\square}}\right)^2\vspace{-.8mm}$ does outperform $\left(\tenq{T}^{(n)}_{\text{\rm W}{\pm}}\right)^2$.\vspace{-4mm}

\subsubsection{Nonspherical Gaussian samples}\label{sec415}\vspace{-3mm}

\begin{figure}[h]
\centering
\begin{tikzpicture}
    \begin{groupplot}[group style={
                      group name=myplot,
                      group size= 3 by 1
                      },height=5cm,width=4.5cm,
                      ytick={0,0.2,0.4,0.6,0.8,1},
                      ymax=1.1]
\nextgroupplot[title={$n=100$}]\addplot [domain=0:.5, color=blue] coordinates { 
(0,0.042)
(0.1,0.068)
(0.2,0.128)
(0.3,0.216)
(0.4,0.424)
(0.5,0.636)
}; 
\label{plots:plot1} 
\addplot [domain=0:.5, color=red] coordinates { 
(0,0.046)
(0.1,0.058)
(0.2,0.124)
(0.3,0.222)
(0.4,0.434)
(0.5,0.638)
}; 
\label{plots:plot2} 
\addplot [domain=0:.5, color=orange] coordinates { 
(0,0.07)
(0.1,0.056)
(0.2,0.122)
(0.3,0.244)
(0.4,0.434)
(0.5,0.628)
}; 
\label{plots:plotH} 
\addplot [domain=0:.5, color=grey] coordinates { 
(0,0.05)
(0.1,0.05)
(0.2,0.05)
(0.3,0.05)
(0.4,0.05)
(0.5,0.05)
}; 
\label{plots:plot3} 
\nextgroupplot[title={$n=200$}]\addplot [domain=0:.5, color=blue] coordinates { 
(0,0.048)
(0.1,0.102)
(0.2,0.26)
(0.3,0.48)
(0.4,0.744000000000001)
(0.5,0.930000000000001)
}; 
\addplot [domain=0:.5, color=red] coordinates { 
(0,0.048)
(0.1,0.0960000000000001)
(0.2,0.252)
(0.3,0.45)
(0.4,0.734000000000001)
(0.5,0.938000000000001)
}; 
\addplot [domain=0:.5, color=orange] coordinates { 
(0,0.05)
(0.1,0.104)
(0.2,0.30)
(0.3,0.486)
(0.4,0.73)
(0.5,0.902)
}; 
\addplot [domain=0:.5, color=grey] coordinates { 
(0,0.05)
(0.1,0.05)
(0.2,0.05)
(0.3,0.05)
(0.4,0.05)
(0.5,0.05)
}; 
\nextgroupplot[title={$n=400$}]\addplot [domain=0:.5, color=blue] coordinates { 
(0,0.038)
(0.1,0.148)
(0.2,0.44)
(0.3,0.808000000000001)
(0.4,0.976000000000001)
(0.5,0.998000000000001)
}; 
\addplot [domain=0:.5, color=red] coordinates { 
(0,0.052)
(0.1,0.158)
(0.2,0.43)
(0.3,0.800000000000001)
(0.4,0.966000000000001)
(0.5,0.998000000000001)
}; 
\addplot [domain=0:.5, color=orange] coordinates { 
(0,0.064)
(0.1,0.146)
(0.2,0.482)
(0.3,0.826)
(0.4,0.978)
(0.5,0.998)
}; 
\addplot [domain=0:.5, color=grey] coordinates { 
(0,0.05)
(0.1,0.05)
(0.2,0.05)
(0.3,0.05)
(0.4,0.05)
(0.5,0.05)
}; 

            \end{groupplot}      
\path (myplot c1r1.south west|-current bounding box.south)--
      coordinate(legendpos)
      (myplot c3r1.south east|-current bounding box.south);
\matrix[
    matrix of nodes,
    anchor=south,
    draw,
    inner sep=0.2em,
    draw
  ]at([yshift=-6ex]legendpos)
  {
    \ref{plots:plot1}& $\left(\tenq{T}^{(n)}_{\text{\rm W}{\pm}}\right)^2$ &[5pt]
    \ref{plots:plot2}& $\left(\tenq{T}^{(n)}_{\text{\rm W}{\scriptscriptstyle\square}}\right)^2$ &[5pt]
    \ref{plots:plotH}& $T^2$
    \\};
\end{tikzpicture}
\caption{\small\slshape Rejection frequencies, for samples with  nonspherical Gaussian distributions (see 4.1.5) and various sample sizes, of Hotelling's test based on $T^2$ and the Wilcoxon tests based on $\left(\tenq{T}^{(n)}_{\text{\rm W}{\pm}}\right)^2\vspace{-.5mm}$ and~$\left(\tenq{T}^{(n)}_{\text{\rm W}{\scriptscriptstyle\square}}\right)^2$, respectively,  as  functions of the shift $\eta$;  
$N=500$ replications.\vspace{-2mm}} 
\label{fig: BivCompGausCorWilcox}
\end{figure}
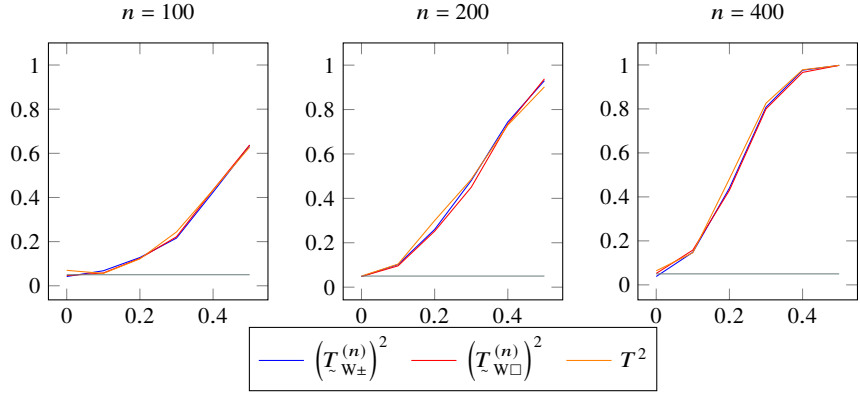

The first sample is drawn from a $ \mathcal{N}( (0,0)\pr, \Sigmab)$ distribution,    the second sample is drawn from a $\mathcal{N}( (\eta,\eta)\pr, \Sigmab)$ one;    $\text{vech}(\Sigmab) = (1,0.8,1)\pr$.   Rejection frequencies over~$N=500$ replications are
 shown (as  functions of $\eta$) in Figure~\ref{fig: BivCompGausCorWilcox}.  The results are essentially the same as in the spherical case (Section~\ref{sec411}). Note the loss of power  in the three tests under study, due to the non-specification of the population covariance matrix; that loss, however, is uniform over the three tests. 


\subsubsection{Spherical Cauchy samples}\label{sec416}\vspace{-3mm}

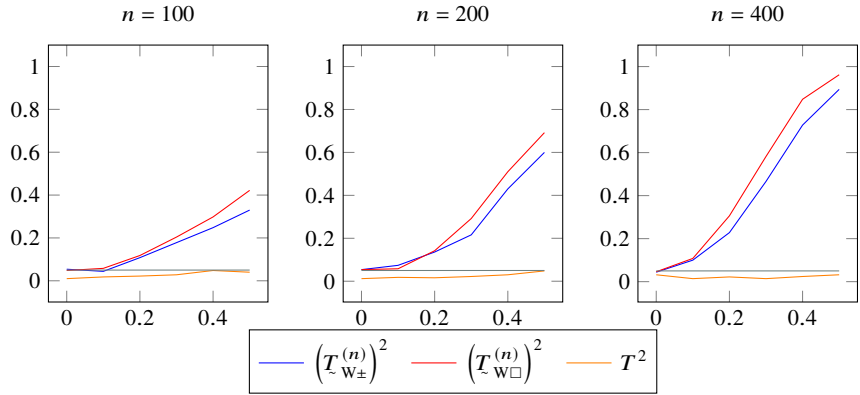
\begin{figure}[h]
\centering
\begin{tikzpicture}
    \begin{groupplot}[group style={
                      group name=myplot,
                      group size= 3 by 1
                      },height=5cm,width=4.5cm,
                      ytick={0,0.2,0.4,0.6,0.8,1},
                      ymax=1.1]
\nextgroupplot[title={$n=100$}]\addplot [domain=0:.5, color=blue] coordinates { 
(0,0.054)
(0.1,0.044)
(0.2,0.108)
(0.3,0.178)
(0.4,0.248)
(0.5,0.33)
}; 
\label{plots:plot1} 
\addplot [domain=0:.5, color=red] coordinates { 
(0,0.048)
(0.1,0.058)
(0.2,0.118)
(0.3,0.204)
(0.4,0.298)
(0.5,0.422)
}; 
\label{plots:plot2} 
\addplot [domain=0:.5, color=orange] coordinates { 
(0,0.01)
(0.1,0.018)
(0.2,0.022)
(0.3,0.028)
(0.4,0.048)
(0.5,0.04)
}; 
\label{plots:plotH} 
\addplot [domain=0:.5, color=grey] coordinates { 
(0,0.05)
(0.1,0.05)
(0.2,0.05)
(0.3,0.05)
(0.4,0.05)
(0.5,0.05)
}; 
\label{plots:plot3} 
\nextgroupplot[title={$n=200$}]\addplot [domain=0:.5, color=blue] coordinates { 
(0,0.054)
(0.1,0.0740000000000001)
(0.2,0.136)
(0.3,0.216)
(0.4,0.43)
(0.5,0.6)
}; 
\addplot [domain=0:.5, color=red] coordinates { 
(0,0.054)
(0.1,0.058)
(0.2,0.142)
(0.3,0.292)
(0.4,0.51)
(0.5,0.692000000000001)
}; 
\addplot [domain=0:.5, color=orange] coordinates { 
(0,0.012)
(0.1,0.018)
(0.2,0.016)
(0.3,0.022)
(0.4,0.03)
(0.5,0.048)
}; 
\addplot [domain=0:.5, color=grey] coordinates { 
(0,0.05)
(0.1,0.05)
(0.2,0.05)
(0.3,0.05)
(0.4,0.05)
(0.5,0.05)
}; 
\nextgroupplot[title={$n=400$}]\addplot [domain=0:.5, color=blue] coordinates { 
(0,0.044)
(0.1,0.1)
(0.2,0.228)
(0.3,0.466)
(0.4,0.728000000000001)
(0.5,0.894000000000001)
}; 
\addplot [domain=0:.5, color=red] coordinates { 
(0,0.046)
(0.1,0.108)
(0.2,0.306)
(0.3,0.582)
(0.4,0.848000000000001)
(0.5,0.962000000000001)
}; 
\addplot [domain=0:.5, color=orange] coordinates { 
(0,0.032)
(0.1,0.014)
(0.2,0.022)
(0.3,0.014)
(0.4,0.02400000001)
(0.5,0.032000000000001)
}; 
\addplot [domain=0:.5, color=grey] coordinates { 
(0,0.05)
(0.1,0.05)
(0.2,0.05)
(0.3,0.05)
(0.4,0.05)
(0.5,0.05)
}; 
            \end{groupplot}      
\path (myplot c1r1.south west|-current bounding box.south)--
      coordinate(legendpos)
      (myplot c3r1.south east|-current bounding box.south);
\matrix[
    matrix of nodes,
    anchor=south,
    draw,
    inner sep=0.2em,
    draw
  ]at([yshift=-6ex]legendpos)
  {
    \ref{plots:plot1}& $\left(\tenq{T}^{(n)}_{\text{\rm W}{\pm}}\right)^2$ &[5pt]
    \ref{plots:plot2}& $\left(\tenq{T}^{(n)}_{\text{\rm W}{\scriptscriptstyle\square}}\right)^2$&[5pt]
    \ref{plots:plotH}& $T^2$
     \\};
\end{tikzpicture}
\caption{\small\slshape Rejection frequencies, for samples with  spherical Cauchy distributions (see 4.1.6) and various sample sizes, of Hotelling's test based on $T^2$ and  the Wilcoxon tests based on $\left(\tenq{T}^{(n)}_{\text{\rm W}{\pm}}\right)^2\vspace{-.5mm}$ and~$\left(\tenq{T}^{(n)}_{\text{\rm W}{\scriptscriptstyle\square}}\right)^2$, respectively,  as  functions of the shift $\eta$;  
$N=500$ replications.\vspace{-3.5mm}} 
\label{fig: BivCompCsphWilcox}
\end{figure}

The first sample is drawn from a centered spherical Student with one degree of freedom $ t_{1} ( (0,0)\pr, {\bf I}_2)$ (spherical Cauchy),   the second one from the shifted version~$ t_{1} ((\eta, \eta)\pr , {\bf I}_2)$ of the same distribution.   Rejection frequencies over $N=500$ replications are
 shown (as  functions of $\eta$) in Figure~\ref{fig: BivCompCsphWilcox}.  The performance of Hotelling, again, is a disaster; although the actual distribution is spherical, $\left(\tenq{T}^{(n)}_{\text{\rm W}{\scriptscriptstyle\square}}\right)^2\vspace{-.8mm}$ still outperforms $\left(\tenq{T}^{(n)}_{\text{\rm W}{\pm}}\right)^2$. \vspace{-3mm}

\subsection{Wilcoxon-type statistics in dimension $d=5$} \vspace{-2.5mm}

We essentially adopted the same simulation settings as before, with $n_1=n_2=n/2$. A sample size of $n=100$ in dimension $d=5$ is very small, though, and we considered sample   sizes~$n=200$, 400, and~800. \vspace{-3mm}




\subsubsection{Spherical Gaussian samples}\label{sec421}\vspace{-3mm}

\begin{figure}[h]
\centering
\begin{tikzpicture}
    \begin{groupplot}[group style={
                      group name=myplot,
                      group size= 3 by 1
                      },height=5cm,width=4.5cm,
                      ytick={0,0.2,0.4,0.6,0.8,1},
                      ymax=1.1]
\nextgroupplot[title={$n=200$}]\addplot [domain=0:.5, color=blue] coordinates { 
(0,0.046)
(0.1,0.138)
(0.2,0.444)
(0.3,0.878000000000001)
(0.4,0.994000000000001)
(0.5,1)
}; 
\addplot [domain=0:.5, color=red] coordinates { 
(0,0.052)
(0.1,0.136)
(0.2,0.458)
(0.3,0.862000000000001)
(0.4,0.996000000000001)
(0.5,1)
}; 
\addplot [domain=0:.5, color=orange] coordinates { 
(0,0.064)
(0.1,0.188)
(0.2,0.676)
(0.3,0.966000000000001)
(0.4,1)
(0.5,1)
}; 
\addplot [domain=0:.5, color=grey] coordinates { 
(0,0.05)
(0.1,0.05)
(0.2,0.05)
(0.3,0.05)
(0.4,0.05)
(0.5,0.05)
}; 
\nextgroupplot[title={$n=400$}]\addplot [domain=0:.5, color=blue] coordinates { 
(0,0.052)
(0.1,0.246)
(0.2,0.860000000000001)
(0.3,1)
(0.4,1)
(0.5,1)
}; 
\addplot [domain=0:.5, color=red] coordinates { 
(0,0.044)
(0.1,0.29)
(0.2,0.884000000000001)
(0.3,1)
(0.4,1)
(0.5,1)
}; 
\addplot [domain=0:.5, color=orange] coordinates { 
(0,0.056)
(0.1,0.4)
(0.2,0.952)
(0.3,1)
(0.4,1)
(0.5,1)
}; 
\addplot [domain=0:.5, color=grey] coordinates { 
(0,0.05)
(0.1,0.05)
(0.2,0.05)
(0.3,0.05)
(0.4,0.05)
(0.5,0.05)
}; 
\nextgroupplot[title={$n=800$}]\addplot [domain=0:.5, color=blue] coordinates { 
(0,0.042)
(0.05,0.158)
(0.1,0.588)
(0.15,0.926000000000001)
(0.2,0.996000000000001)
(0.25,1)
(0.5,1)
}; 
\label{plots:plot1Wil5D5} 
\addplot [domain=0:.5, color=red] coordinates { 
(0,0.058)
(0.05,0.166)
(0.1,0.59)
(0.15,0.916000000000001)
(0.2,0.998000000000001)
(0.25,1)
(0.5,1)
}; 
\label{plots:plot2Wil5D5} 
\addplot [domain=0:.5, color=orange] coordinates { 
(0,0.038)
(0.05,0.194)
(0.1,0.68)
(0.15,0.968000000000001)
(0.2,0.998000000000001)
(0.25,1)
(0.5,1)
}; 
\label{plots:plot3Wil5D5} 
\addplot [domain=0:.5, color=grey] coordinates { 
(0,0.05)
(0.05,0.05)
(0.1,0.05)
(0.15,0.05)
(0.2,0.05)
(0.25,0.05)
(0.5,0.05)
}; 
\label{plots:plot4Wil5D5}

            \end{groupplot}      
\path (myplot c1r1.south west|-current bounding box.south)--
      coordinate(legendpos)
      (myplot c3r1.south east|-current bounding box.south);
\matrix[
    matrix of nodes,
    anchor=south,
    draw,
    inner sep=0.2em,
    draw
  ]at([yshift=-6ex]legendpos)
  {
 \ref{plots:plot1}& $\left(\tenq{T}^{(n)}_{\text{\rm W}{\pm}}\right)^2$ &[5pt]
    \ref{plots:plot2}& $\left(\tenq{T}^{(n)}_{\text{\rm W}{\scriptscriptstyle\square}}\right)^2$ &[5pt]
     \ref{plots:plotH}& $T^2$
    \\};
    \end{tikzpicture}
\caption{\small\slshape Rejection frequencies, for samples with 5-dimensional spherical Gaussian distributions (see~4.2.1) and various sample sizes, of Hotelling's test based on $T^2$ and the Wilcoxon tests based on $\left(\tenq{T}^{(n)}_{\text{\rm W}{\pm}}\right)^2$ and~$\left(\tenq{T}^{(n)}_{\text{\rm W}{\scriptscriptstyle\square}}\right)^2$, respectively,  as  functions of the shift $\eta$;   $N=500$ replications.\vspace{-10mm}}
\label{fig: BivCompGausWilcox5D}
\end{figure}
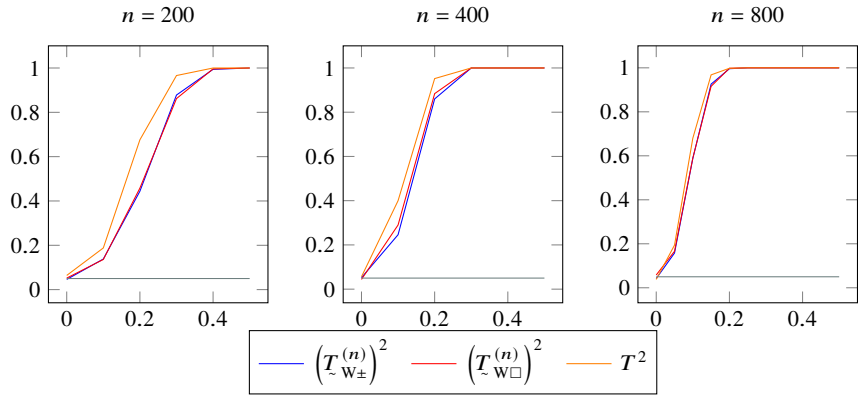

Here, the first sample is drawn from the $ \mathcal{N} ( {\bf 0}, {\bf I}_5)$ distribution, the second one   from  the~$ \mathcal{N} ( {\eta}{\mathbf 1} , {\bf I}_5)$ distribution, where~${\mathbf 1}$ denotes a 5-variate vector of ones.   Rejection frequencies over~$N=500$ replications are
 shown (as  functions of $\eta$) in  Figure~\ref{fig: BivCompGausWilcox5D}.  In the ``small sample'' case $(n=200)$, the optimality of  Hotelling over Wilcoxon is perceptible (more so than in dimension $d=2$); this superiority, however, fades away with growing $n$: again, under Gaussian assumptions, abandoning the parametrically optimal Hotelling test in favor of the rank-based Wilcoxon one has no visible cost in terms of power.    \vspace{-4mm}

\subsubsection{Spherical Student samples}\label{sec422}\vspace{-3mm}

The first sample is drawn from the centered spherical Student with 2.1 degrees of freedom $ t_{2.1} ({\boldsymbol 0}, {\bf I}_5)$, the second one   from  the shifted~$ t_{2.1} ({\eta}{\mathbf 1} , {\bf I}_5)$ distribution, where~${\mathbf 1}$ denotes a 5-variate vector of ones.   Rejection frequencies over~$N=500$ replications 
\begin{figure}[h]
\centering
\begin{tikzpicture}
    \begin{groupplot}[group style={
                      group name=myplot,
                      group size= 3 by 1
                      },height=5cm,width=4.5cm,
                      ytick={0,0.2,0.4,0.6,0.8,1},
                      ymax=1.1]
\nextgroupplot[title={$n=200$}]\addplot [domain=0:.5, color=blue] coordinates { 
(0,0.05)
(0.1,0.0860000000000001)
(0.2,0.274)
(0.3,0.604)
(0.4,0.854000000000001)
(0.5,0.972000000000001)
}; 
\addplot [domain=0:.5, color=red] coordinates { 
(0,0.042)
(0.1,0.11)
(0.2,0.356)
(0.3,0.702000000000001)
(0.4,0.926000000000001)
(0.5,0.988000000000001)
}; 
\addplot [domain=0:.5, color=orange] coordinates { 
(0,0.05)
(0.1,0.068)
(0.2,0.176)
(0.3,0.422)
(0.4,0.616)
(0.5,0.78)
}; 
\addplot [domain=0:.5, color=grey] coordinates { 
(0,0.05)
(0.1,0.05)
(0.2,0.05)
(0.3,0.05)
(0.4,0.05)
(0.5,0.05)
}; 
\nextgroupplot[title={$n=400$}]\addplot [domain=0:.5, color=blue] coordinates { 
(0,0.042)
(0.1,0.132)
(0.2,0.488)
(0.3,0.910000000000001)
(0.4,0.990000000000001)
(0.5,1)
}; 
\addplot [domain=0:.5, color=red] coordinates { 
(0,0.06)
(0.1,0.178)
(0.2,0.66)
(0.3,0.976000000000001)
(0.4,1)
(0.5,1)
}; 
\addplot [domain=0:.5, color=orange] coordinates { 
(0,0.034)
(0.1,0.11)
(0.2,0.294)
(0.3,0.614)
(0.4,0.876)
(0.5,0.95)
}; 
\addplot [domain=0:.5, color=grey] coordinates { 
(0,0.05)
(0.1,0.05)
(0.2,0.05)
(0.3,0.05)
(0.4,0.05)
(0.5,0.05)
}; 
\nextgroupplot[title={$n=800$}]\addplot [domain=0:.5, color=blue] coordinates { 
(0,0.03)
(0.05,0.0920000000000001)
(0.1,0.316)
(0.15,0.588)
(0.2,0.916000000000001)
(0.25,0.986000000000001)
(0.3, 0.998)
(0.4,1)
(0.5,1)
}; 
\label{plots:plot1Wil5D6} 
\addplot [domain=0:.5, color=red] coordinates { 
(0,0.048)
(0.05,0.114)
(0.1,0.378)
(0.15,0.746)
(0.2,0.970000000000001)
(0.25,0.998000000000001)
(0.3, 1)
(0.4,1)
(0.5,1)
}; 
\label{plots:plot2Wil5D6} 
\addplot [domain=0:.5, color=orange] coordinates { 
(0,0.034)
(0.05,0.064)
(0.1,0.124)
(0.15,0.29)
(0.2,0.532)
(0.25,0.734)
(0.3, 0.886)
(0.4,0.948)
(0.5,0.994)
}; 
\label{plots:plot3Wil5D6} 
\addplot [domain=0:.5, color=grey] coordinates { 
(0,0.05)
(0.05,0.05)
(0.1,0.05)
(0.15,0.05)
(0.2,0.05)
(0.25,0.05)
(0.5,0.05)
}; 
\label{plots:plot4Wil5D6} 
            \end{groupplot}      
\path (myplot c1r1.south west|-current bounding box.south)--
      coordinate(legendpos)
      (myplot c3r1.south east|-current bounding box.south);
\matrix[
    matrix of nodes,
    anchor=south,
    draw,
    inner sep=0.2em,
    draw
  ]at([yshift=-6ex]legendpos)
 {
    \ref{plots:plot1}& $\left(\tenq{T}^{(n)}_{\text{\rm W}{\pm}}\right)^2$ &[5pt]
    \ref{plots:plot2}& $\left(\tenq{T}^{(n)}_{\text{\rm W}{\scriptscriptstyle\square}}\right)^2$ &[5pt]
    \ref{plots:plotH}& $T^2$
     \\};
    \end{tikzpicture}
    \caption{\small\slshape Rejection frequencies, for samples with 5-dimensional spherical Student $t_{2.1}$  distributions (see~4.2.2) and various sample sizes, of the Wilcoxon tests based on $\left(\tenq{T}^{(n)}_{\text{\rm W}{\pm}}\right)^2$ and $\left(\tenq{T}^{(n)}_{\text{\rm W}{\scriptscriptstyle\square}}\right)^2$, respectively,  as  functions of the shift $\eta$;   $N=500$ replications.\vspace{-2mm}}
\label{fig: CompTsphWilcox5D}
\end{figure}
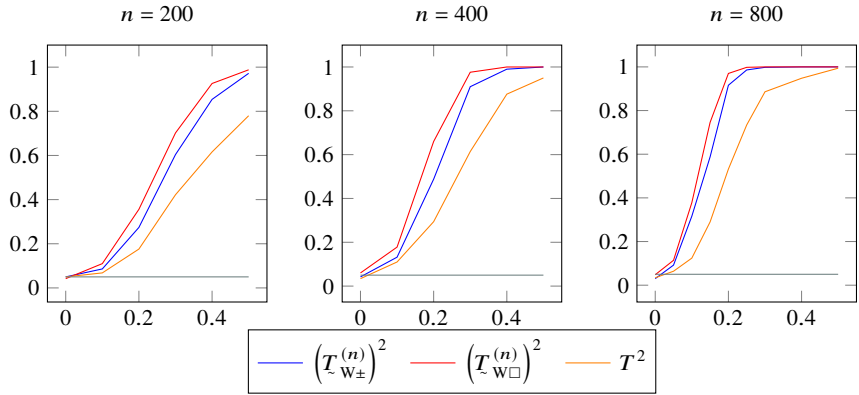
are
 shown (as  functions of $\eta$) in  Figure~\ref{fig: CompTsphWilcox5D}.  The conclusions are quite similar to those
 in dimension $d=2$: although the moments of order two still are finite, the power of Hotelling deteriorates with respect to the Gaussian case. Despite the spherical nature of the actual distribution, slight advantage of $\left(\tenq{T}^{(n)}_{\text{\rm W}{\scriptscriptstyle\square}}\right)^2\vspace{-.8mm}$ over $\left(\tenq{T}^{(n)}_{\text{\rm W}{\pm}}\right)^2$.  \vspace{-6mm}


\subsubsection{Nonspherical Gaussian samples}\label{sec423}\vspace{-3mm}

The first sample is drawn from the $ \mathcal{N} ( {\bf 0}, \Sigmab)$ distribution, the second one   from  the~$ \mathcal{N} ( {\eta}{\mathbf 1} , \Sigmab)$ distribution, where~${\mathbf 1}$ denotes a 5-variate vector of ones and  $\Sigmab$ is a correlation matrix with all off-diagonal entries equal to 0.5.    Rejection frequencies over~$N=500$ replications are
 shown (as  functions of $\eta$) in  Figure~\ref{fig: CompGausCorrWilcox5D}.  Here again, the slight advantage of Hotelling over Wilcoxon very rapidly fades away with growing~$n$, and the three tests yield very similar performances; in parti\-cular, no significant difference can be detected between   $\left(\tenq{T}^{(n)}_{\text{\rm W}{\pm}}\right)^2$ and $\left(\tenq{T}^{(n)}_{\text{\rm W}{\scriptscriptstyle\square}}\right)^2\vspace{-.8mm}$. \vspace{-4mm}

\begin{figure}[h]
\centering
\begin{tikzpicture}
    \begin{groupplot}[group style={
                      group name=myplot,
                      group size= 3 by 1
                      },height=5cm,width=4.5cm,
                      ytick={0,0.2,0.4,0.6,0.8,1},
                      ymax=1.1]
\nextgroupplot[title={$n=200$}]\addplot [domain=0:.5, color=blue] coordinates { 
(0,0.028)
(0.1,0.07)
(0.2,0.16)
(0.3,0.428)
(0.4,0.724000000000001)
(0.5,0.906000000000001)
}; 
\addplot [domain=0:.5, color=red] coordinates { 
(0,0.038)
(0.1,0.0780000000000001)
(0.2,0.188)
(0.3,0.414)
(0.4,0.720000000000001)
(0.5,0.924000000000001)
}; 
\addplot [domain=0:.5, color=orange] coordinates { 
(0,0.042)
(0.1,0.0780000000000001)
(0.2,0.218)
(0.3,0.52)
(0.4,0.81)
(0.5,0.958)
}; 
\addplot [domain=0:.5, color=grey] coordinates { 
(0,0.05)
(0.1,0.05)
(0.2,0.05)
(0.3,0.05)
(0.4,0.05)
(0.5,0.05)
}; 
\nextgroupplot[title={$n=400$}]\addplot [domain=0:.5, color=blue] coordinates { 
(0,0.06)
(0.1,0.128)
(0.2,0.358)
(0.3,0.802000000000001)
(0.4,0.958000000000001)
(0.5,1)
}; 
\addplot [domain=0:.5, color=red] coordinates { 
(0,0.058)
(0.1,0.118)
(0.2,0.422)
(0.3,0.842000000000001)
(0.4,0.970000000000001)
(0.5,1)
}; 
\addplot [domain=0:.5, color=orange] coordinates { 
(0,0.058)
(0.1,0.126)
(0.2,0.488)
(0.3,0.852000000000001)
(0.4,0.996)
(0.5,1)
}; 
\addplot [domain=0:.5, color=grey] coordinates { 
(0,0.05)
(0.1,0.05)
(0.2,0.05)
(0.3,0.05)
(0.4,0.05)
(0.5,0.05)
}; 
\nextgroupplot[title={$n=800$}]\addplot [domain=0:.5, color=blue] coordinates { 
(0,0.058)
(0.05,0.076)
(0.1,0.232)
(0.15,0.504)
(0.2,0.748)
(0.25,0.948000000000001)
(0.3,0.986)
(0.4,1)
(0.5,1)
}; 
\label{plots:plot1Wil5D7} 
\addplot [domain=0:.5, color=red] coordinates { 
(0,0.062)
(0.05,0.0940000000000001)
(0.1,0.254)
(0.15,0.518)
(0.2,0.778000000000001)
(0.25,0.958000000000001)
(0.3,0.992)
(0.4,1)
(0.5,1)
}; 
\label{plots:plot2Wil5D7} 
\addplot [domain=0:.5, color=orange] coordinates { 
(0,0.044)
(0.05,0.084)
(0.1,0.264)
(0.15,0.498)
(0.2,0.824)
(0.25,0.954000000000001)
(0.3,1)
(0.5,1)
}; 
\label{plots:plot3Wil5D7} 
\addplot [domain=0:.5, color=grey] coordinates { 
(0,0.05)
(0.05,0.05)
(0.1,0.05)
(0.15,0.05)
(0.2,0.05)
(0.25,0.05)
}; 
\label{plots:plot4Wil5D7}

            \end{groupplot}      
\path (myplot c1r1.south west|-current bounding box.south)--
      coordinate(legendpos)
      (myplot c3r1.south east|-current bounding box.south);
\matrix[
    matrix of nodes,
    anchor=south,
    draw,
    inner sep=0.2em,
    draw
  ]at([yshift=-6ex]legendpos)
 {
 \ref{plots:plot1}& $\left(\tenq{T}^{(n)}_{\text{\rm W}{\pm}}\right)^2$ &[5pt]
    \ref{plots:plot2}& $\left(\tenq{T}^{(n)}_{\text{\rm W}{\scriptscriptstyle\square}}\right)^2$ &[5pt]
     \ref{plots:plotH}& $T^2$
     \\};
\end{tikzpicture}
  \caption{\small\slshape Rejection frequencies, for nonspherical 5-dimensional Gaussian  distributions (see~4.2.3) and various sample sizes, of the Wilcoxon tests based on $\left(\tenq{T}^{(n)}_{\text{\rm W}{\pm}}\right)^2$ and $\left(\tenq{T}^{(n)}_{\text{\rm W}{\scriptscriptstyle\square}}\right)^2$, respectively,  as  functions of the shift $\eta$;   $N=500$ replications.\vspace{-8mm}}
\label{fig: CompGausCorrWilcox5D}
\end{figure}
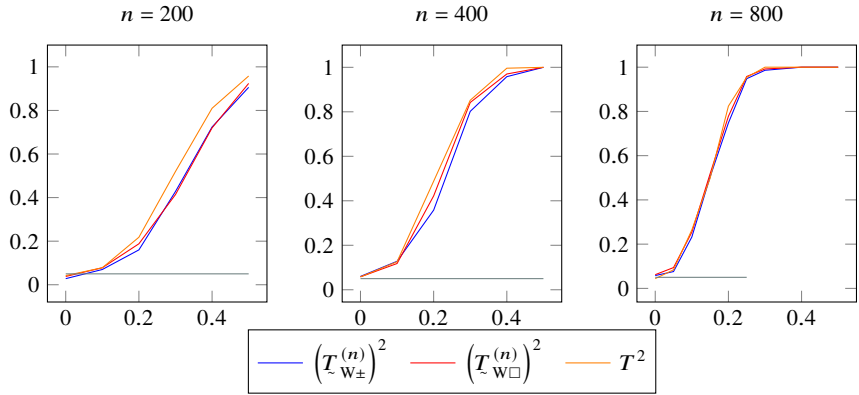
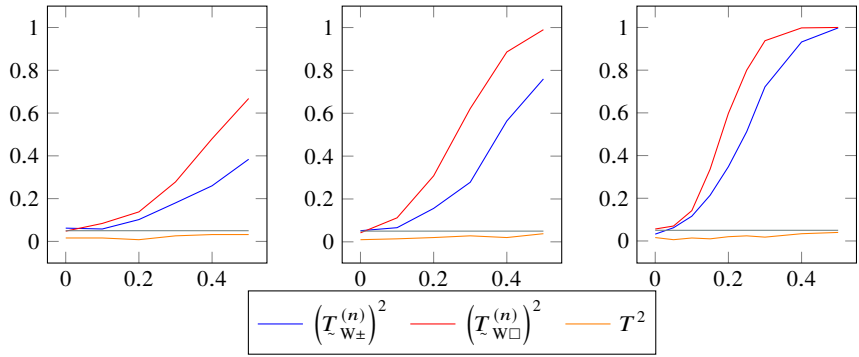
\begin{figure}[h]
\centering
\begin{tikzpicture}
    \begin{groupplot}[group style={
                      group name=myplot,
                      group size= 3 by 1
                      },height=5cm,width=4.5cm,
                      ytick={0,0.2,0.4,0.6,0.8,1},
                      ymax=1.1]
\nextgroupplot[title={$n=200$}]\addplot [domain=0:.5, color=blue] coordinates { 
(0,0.062)
(0.1,0.058)
(0.2,0.102)
(0.3,0.18)
(0.4,0.26)
(0.5,0.384)
}; 
\addplot [domain=0:.5, color=red] coordinates { 
(0,0.048)
(0.1,0.0840000000000001)
(0.2,0.138)
(0.3,0.278)
(0.4,0.48)
(0.5,0.668)
}; 
\addplot [domain=0:.5, color=orange] coordinates { 
(0,0.016)
(0.1,0.016)
(0.2,0.008)
(0.3,0.026)
(0.4,0.032)
(0.5,0.032)
}; 
\addplot [domain=0:.5, color=grey] coordinates { 
(0,0.05)
(0.1,0.05)
(0.2,0.05)
(0.3,0.05)
(0.4,0.05)
(0.5,0.05)
}; 
\nextgroupplot[title={$n=400$}]\addplot [domain=0:.5, color=blue] coordinates { 
(0,0.052)
(0.1,0.066)
(0.2,0.156)
(0.3,0.278)
(0.4,0.564)
(0.5,0.760000000000001)
}; 
\addplot [domain=0:.5, color=red] coordinates { 
(0,0.042)
(0.1,0.112)
(0.2,0.308)
(0.3,0.622)
(0.4,0.886000000000001)
(0.5,0.990000000000001)
}; 
\addplot [domain=0:.5, color=orange] coordinates { 
(0,0.01)
(0.1,0.014)
(0.2,0.02)
(0.3,0.028)
(0.4,0.02)
(0.5,0.038)
}; 
\addplot [domain=0:.5, color=grey] coordinates { 
(0,0.05)
(0.1,0.05)
(0.2,0.05)
(0.3,0.05)
(0.4,0.05)
(0.5,0.05)
}; 
\nextgroupplot[title={$n=800$}]\addplot [domain=0:.5, color=blue] coordinates { 
(0,0.032)
(0.05,0.062)
(0.1,0.116)
(0.15,0.214)
(0.2,0.348)
(0.25,0.512)
(0.3,0.722)
(0.4,0.932)
(0.5,0.998)
}; 
\label{plots:plot1Wil5D8} 
\addplot [domain=0:.5, color=red] coordinates { 
(0,0.056)
(0.05,0.07)
(0.1,0.142)
(0.15,0.336)
(0.2,0.598)
(0.25,0.800000000000001)
(0.3,0.938)
(0.4,0.998)
(0.5,1)
}; 
\label{plots:plot2Wil5D8} 
\addplot [domain=0:.5, color=orange] coordinates { 
(0,0.016)
(0.05,0.006)
(0.1,0.014)
(0.15,0.01)
(0.2,0.02)
(0.25,0.024)
(0.3,0.018)
(0.4,0.034)
(0.5,0.04)
}; 
\label{plots:plot3Wil5D8} 
\addplot [domain=0:.5, color=grey] coordinates { 
(0,0.05)
(0.05,0.05)
(0.1,0.05)
(0.15,0.05)
(0.2,0.05)
(0.25,0.05)
(0.3,0.05)
(0.4,0.05)
(0.5,0.05)
}; 
\label{plots:plot4Wil5D8}

            \end{groupplot}      
\path (myplot c1r1.south west|-current bounding box.south)--
      coordinate(legendpos)
      (myplot c3r1.south east|-current bounding box.south);
\matrix[
    matrix of nodes,
    anchor=south,
    draw,
    inner sep=0.2em,
    draw
  ]at([yshift=-6ex]legendpos)
 {
 \ref{plots:plot1}& $\left(\tenq{T}^{(n)}_{\text{\rm W}{\pm}}\right)^2$ &[5pt]
    \ref{plots:plot2}& $\left(\tenq{T}^{(n)}_{\text{\rm W}{\scriptscriptstyle\square}}\right)^2$&[5pt]
     \ref{plots:plotH}& $T^2$
     \\};
\end{tikzpicture}
  \caption{\small\slshape Rejection frequencies, for 5-dimensional   distributions with independent Cauchy marginals  (see~4.2.4) and various sample sizes, of Hotelling's test based on $T^2$ and the Wilcoxon tests based on $\left(\tenq{T}^{(n)}_{\text{\rm W}{\pm}}\right)^2\vspace{-.6mm}$ and~$\left(\tenq{T}^{(n)}_{\text{\rm W}{\scriptscriptstyle\square}}\right)^2$, respectively,  as  functions of the shift $\eta$;   $N=500$ replications.\vspace{-3mm}}\label{fig: CompCauchyWilcox5D}
\end{figure}

\vspace{-0mm}



\subsubsection{Samples with independent Cauchy marginals}\label{sec424}\vspace{-3mm}

The first sample is drawn from a product of five independent Cauchy distributions,   the second one from the shifted version of the same.   Rejection frequencies over~$500$ replications are
 shown (as  functions of $\eta$) in  Figure~\ref{fig: CompCauchyWilcox5D}.  The performance of Hotelling, as in dimension $d=2$, is terrible. The advantage (which is in line with the independent component nature of the distribution) of $\left(\tenq{T}^{(n)}_{\text{\rm W}{\scriptscriptstyle\square}}\right)^2\vspace{-.8mm}$ over  $\left(\tenq{T}^{(n)}_{\text{\rm W}{\pm}}\right)^2\vspace{1mm}$ is even more significant than in dimension $d=2$.\vspace{-5mm}

%
%

\section{van der Waerden-type  tests}\vspace{-1mm}
\label{sec: Scores}

Below we are considering four distinct tests of the van der Waerden type, based (see Section~\ref{sec:teststat}) on $\left(\tenq{T}^{(n)}_{\text{\rm vdW}{\pm}}\right)^2\vspace{-1mm}$ (spherical uniform reference density, ($\mathfrak G$i) grid),~$\left(\tenq{T}^{(n)}_{\text{\rm vdW}{\scriptscriptstyle\square}}\right)^2$ (cubic uniform reference density, ($\mathfrak G$ii) grid), $\left(\tenq{T}^{(n)}_{\text{\rm vdW}{\pm}{\mathcal N}}\right)^2$ (spherical Gaussian reference density, spherical grid ($\mathfrak G$iii)), and $\left(\tenq{T}^{(n)}_{\text{\rm vdW}{\scriptscriptstyle\square}{\mathcal N}}\right)^2\vspace{.6mm}$  (spherical Gaussian reference density, cubic grid ($\mathfrak G$iv)). \vspace{-4mm}
\subsection{Bivariate case}
\label{sec: VdW2D}

%
%
%

\subsubsection{Spherical Gaussian samples}\label{sec511}\vspace{-3mm}

Same Gaussian samples as in Section~\ref{sec411}. 
Rejection frequencies over~$N=500$ replications are
 shown (as  functions of $\eta$) in  Figure~\ref{fig: BivCompGausVDW}. Indistinctiveness between the performances of Hotelling and the various rank-based tests is even more pronounced than for the Wilcoxon tests: definitely, performing rank-based van der Waerden tests does not imply any loss of efficiency in the Gaussian case.\vspace{-4mm}

\begin{figure}[h]
\centering
\begin{tikzpicture}
    \begin{groupplot}[group style={
                      group name=myplot,
                      group size= 3 by 1
                      },height=5cm,width=4.5cm,
                      ytick={0,0.2,0.4,0.6,0.8,1},
                      ymax=1.1]
\nextgroupplot[title={$n=100$}]\addplot [domain=0:.5, color=black] coordinates { 
(0,0.054)
(0.1,0.0820000000000001)
(0.2,0.198)
(0.3,0.458)
(0.4,0.652)
(0.5,0.838000000000001)
}; 
\label{plots:plot1VDW} 
\addplot [domain=0:.5, color=brightgreen] coordinates { 
(0,0.05)
(0.1,0.0880000000000001)
(0.2,0.204)
(0.3,0.454)
(0.4,0.656)
(0.5,0.848000000000001)
}; 
\label{plots:plot2VDW} 
\addplot [domain=0:.5, color=blue] coordinates { 
(0,0.052)
(0.1,0.0880000000000001)
(0.2,0.196)
(0.3,0.464)
(0.4,0.676)
(0.5,0.848000000000001)
}; 
\label{plots:plot3VDW} 
\addplot [domain=0:.5, color=orange] coordinates { 
(0, 0.052)
(0.1, 0.108)
(0.2, 0.218)
(0.3, 0.472)
(0.4, 0.714)
(0.5, 0.872)
}; 
\label{plots:plotHVDW} 
\addplot [domain=0:.5, color=red] coordinates { 
(0,0.054)
(0.1,0.0900000000000001)
(0.2,0.19)
(0.3,0.452)
(0.4,0.666)
(0.5,0.830000000000001)
}; 
\label{plots:plot4VDW} 
\addplot [domain=0:.5, color=grey] coordinates { 
(0,0.05)
(0.1,0.05)
(0.2,0.05)
(0.3,0.05)
(0.4,0.05)
(0.5,0.05)
}; 
\label{plots:plot5VDW} 
\nextgroupplot[title={$n=200$}]\addplot [domain=0:.5, color=black] coordinates { 
(0,0.054)
(0.1,0.116)
(0.2,0.38)
(0.3,0.752000000000001)
(0.4,0.950000000000001)
(0.5,0.994000000000001)
}; 
\addplot [domain=0:.5, color=brightgreen] coordinates { 
(0,0.06)
(0.1,0.12)
(0.2,0.39)
(0.3,0.756000000000001)
(0.4,0.952000000000001)
(0.5,0.992000000000001)
}; 
\addplot [domain=0:.5, color=blue] coordinates { 
(0,0.05)
(0.1,0.118)
(0.2,0.382)
(0.3,0.744000000000001)
(0.4,0.956000000000001)
(0.5,0.996000000000001)
}; 
\addplot [domain=0:.5, color=orange] coordinates { 
(0, 0.038)
(0.1, 0.114)
(0.2, 0.426)
(0.3, 0.790)
(0.4, .946)
(0.5, 0.994)
}; 

\addplot [domain=0:.5, color=red] coordinates { 
(0,0.046)
(0.1,0.124)
(0.2,0.386)
(0.3,0.746000000000001)
(0.4,0.952000000000001)
(0.5,0.992000000000001)
}; 
\addplot [domain=0:.5, color=grey] coordinates { 
(0,0.05)
(0.1,0.05)
(0.2,0.05)
(0.3,0.05)
(0.4,0.05)
(0.5,0.05)
}; 
\nextgroupplot[title={$n=400$}]\addplot [domain=0:.5, color=black] coordinates { 
(0,0.042)
(0.1,0.192)
(0.2,0.704000000000001)
(0.3,0.960000000000001)
(0.4,1)
(0.5,1)
}; 
\addplot [domain=0:.5, color=brightgreen] coordinates { 
(0,0.046)
(0.1,0.214)
(0.2,0.698000000000001)
(0.3,0.966000000000001)
(0.4,1)
(0.5,1)
}; 
\addplot [domain=0:.5, color=blue] coordinates { 
(0,0.048)
(0.1,0.214)
(0.2,0.702000000000001)
(0.3,0.966000000000001)
(0.4,1)
(0.5,1)
}; 
\addplot [domain=0:.5, color=orange] coordinates { 
(0, 0.066)
(0.1, 0.220)
(0.2, 0.736)
(0.3, 0.974)
(0.4, 1)
(0.5, 1)
}; 
\addplot [domain=0:.5, color=red] coordinates { 
(0,0.042)
(0.1,0.196)
(0.2,0.692000000000001)
(0.3,0.964000000000001)
(0.4,1)
(0.5,1)
}; 
\addplot [domain=0:.5, color=grey] coordinates { 
(0,0.05)
(0.1,0.05)
(0.2,0.05)
(0.3,0.05)
(0.4,0.05)
(0.5,0.05)
}; 

            \end{groupplot}      
\path (myplot c1r1.south west|-current bounding box.south)--
      coordinate(legendpos)
      (myplot c3r1.south east|-current bounding box.south);
\matrix[
    matrix of nodes,
    anchor=south,
    draw,
    inner sep=0.2em,
    draw
  ]at([yshift=-6ex]legendpos)
  {
     \ref{plots:plot3VDW}& $\left(\tenq{T}^{(n)}_{\text{\rm vdW}{\pm}}\right)^2$&[5pt]
        \ref{plots:plot4VDW}&$\left(\tenq{T}^{(n)}_{\text{\rm vdW}{\scriptscriptstyle\square}}\right)^2$ &[5pt]
    \ref{plots:plot2VDW}& $\left(\tenq{T}^{(n)}_{\text{\rm vdW}{\pm}{\mathcal N}}\right)^2$&[5pt]
             \ref{plots:plot1VDW}& $\left(\tenq{T}^{(n)}_{\text{\rm vdW}{\scriptscriptstyle\square}{\mathcal N}}\right)^2$&[5pt]
             \ref{plots:plotHVDW}& $T^2$
    \\};
\end{tikzpicture}
\caption{\small\slshape Rejection frequencies, for  spherical Gaussian samples (see~5.1.1) and various sample sizes, of Hotelling's test based on $T^2$ and the van der Waerden tests based on $\left(\tenq{T}^{(n)}_{\text{\rm vdW}{\pm}}\right)^2$,  $\left(\tenq{T}^{(n)}_{\text{\rm vdW}{\scriptscriptstyle\square}}\right)^2$, $\left(\tenq{T}^{(n)}_{\text{\rm vdW}{\pm}{\mathcal N}}\right)^2$, 
and $\left(\tenq{T}^{(n)}_{\text{\rm vdW}{\scriptscriptstyle\square}{\mathcal N}}\right)^2$, respectively,  as  functions of the shift $\eta$;   $N=500$ replications.\vspace{-2mm}}
\label{fig: BivCompGausVDW}
\end{figure}


\subsubsection{Spherical Student samples}\label{sec512}\vspace{-3mm}

Same Student samples as in Section~\ref{sec412}. 
Rejection frequencies over~$N=500$ replications are
 shown (as  functions of $\eta$) in  Figure~\ref{fig: BivCompTsphVDW}.
   As in dimension $d=2$, powers are dropping (compared with the Gaussian case). The power of Hotelling, however,  deteriorates much more than that   of the various versions of van der Waerden  tests. The latter all yield very similar performance. \vspace{-4mm}

\begin{figure}[h]
\centering
\begin{tikzpicture}
    \begin{groupplot}[group style={
                      group name=myplot,
                      group size= 3 by 1
                      },height=5cm,width=4.5cm,
                      ytick={0,0.2,0.4,0.6,0.8,1},
                      ymax=1.1]
\nextgroupplot[title={$n=100$}]\addplot [domain=0:.5, color=black] coordinates { 
(0,0.058)
(0.1,0.056)
(0.2,0.15)
(0.3,0.236)
(0.4,0.382)
(0.5,0.522)
}; 
\label{plots:plot1VDW} 
\addplot [domain=0:.5, color=brightgreen] coordinates { 
(0,0.056)
(0.1,0.062)
(0.2,0.144)
(0.3,0.24)
(0.4,0.418)
(0.5,0.56)
}; 
\label{plots:plot2VDW} 
\addplot [domain=0:.5, color=blue] coordinates { 
(0,0.058)
(0.1,0.058)
(0.2,0.144)
(0.3,0.244)
(0.4,0.442)
(0.5,0.556)
}; 
\label{plots:plot3VDW} 
\addplot [domain=0:.5, color=orange] coordinates { 
(0, 0.05)
(0.1, 0.062)
(0.2, 0.064)
(0.3, 0.132)
(0.4, 0.228)
(0.5, 0.318)
}; 
\label{plots:plotHVDW} 
\addplot [domain=0:.5, color=red] coordinates { 
(0,0.06)
(0.1,0.054)
(0.2,0.15)
(0.3,0.222)
(0.4,0.398)
(0.5,0.536)
}; 
\label{plots:plot4VDW} 
\addplot [domain=0:.5, color=grey] coordinates { 
(0,0.05)
(0.1,0.05)
(0.2,0.05)
(0.3,0.05)
(0.4,0.05)
(0.5,0.05)
}; 
\label{plots:plot5} 
\nextgroupplot[title={$n=200$}]\addplot [domain=0:.5, color=black] coordinates { 
(0,0.052)
(0.1,0.0860000000000001)
(0.2,0.224)
(0.3,0.462)
(0.4,0.65)
(0.5,0.868000000000001)
}; 
\addplot [domain=0:.5, color=brightgreen] coordinates { 
(0,0.056)
(0.1,0.0880000000000001)
(0.2,0.234)
(0.3,0.484)
(0.4,0.684)
(0.5,0.890000000000001)
}; 
\addplot [domain=0:.5, color=blue] coordinates { 
(0,0.05)
(0.1,0.0980000000000001)
(0.2,0.226)
(0.3,0.478)
(0.4,0.682)
(0.5,0.894000000000001)
}; 
\addplot [domain=0:.5, color=orange] coordinates { 
(0, 0.048)
(0.1, 0.042)
(0.2, 0.09)
(0.3, 0.224)
(0.4, 0.388)
(0.5, 0.496)
}; 
\addplot [domain=0:.5, color=red] coordinates { 
(0,0.054)
(0.1,0.0920000000000001)
(0.2,0.22)
(0.3,0.454)
(0.4,0.656)
(0.5,0.882000000000001)
}; 
\addplot [domain=0:.5, color=grey] coordinates { 
(0,0.05)
(0.1,0.05)
(0.2,0.05)
(0.3,0.05)
(0.4,0.05)
(0.5,0.05)
}; 
\nextgroupplot[title={$n=400$}]\addplot [domain=0:.5, color=black] coordinates { 
(0,0.04)
(0.1,0.112)
(0.2,0.442)
(0.3,0.746000000000001)
(0.4,0.956000000000001)
(0.5,0.992000000000001)
}; 
\addplot [domain=0:.5, color=brightgreen] coordinates { 
(0,0.04)
(0.1,0.122)
(0.2,0.468)
(0.3,0.758000000000001)
(0.4,0.962000000000001)
(0.5,0.994000000000001)
}; 
\addplot [domain=0:.5, color=blue] coordinates { 
(0,0.038)
(0.1,0.126)
(0.2,0.474)
(0.3,0.772000000000001)
(0.4,0.964000000000001)
(0.5,0.994000000000001)
}; 
\addplot [domain=0:.5, color=orange] coordinates { 
(0, 0.028)
(0.1, 0.072)
(0.2, 0.188)
(0.3, 0.352)
(0.4, 0.548)
(0.5, 0.738)
}; 
\addplot [domain=0:.5, color=red] coordinates { 
(0,0.03)
(0.1,0.13)
(0.2,0.462)
(0.3,0.748000000000001)
(0.4,0.954000000000001)
(0.5,0.992000000000001)
}; 
\addplot [domain=0:.5, color=grey] coordinates { 
(0,0.05)
(0.1,0.05)
(0.2,0.05)
(0.3,0.05)
(0.4,0.05)
(0.5,0.05)
}; 
            \end{groupplot}      
\path (myplot c1r1.south west|-current bounding box.south)--
      coordinate(legendpos)
      (myplot c3r1.south east|-current bounding box.south);
\matrix[
    matrix of nodes,
    anchor=south,
    draw,
    inner sep=0.2em,
    draw
  ]at([yshift=-6ex]legendpos)
  {
      \ref{plots:plot3VDW}& $\left(\tenq{T}^{(n)}_{\text{\rm vdW}{\pm}}\right)^2$&[5pt]
        \ref{plots:plot4VDW}&$\left(\tenq{T}^{(n)}_{\text{\rm vdW}{\scriptscriptstyle\square}}\right)^2$ &[5pt]
    \ref{plots:plot2VDW}& $\left(\tenq{T}^{(n)}_{\text{\rm vdW}{\pm}{\mathcal N}}\right)^2$&[5pt]
             \ref{plots:plot1VDW}& $\left(\tenq{T}^{(n)}_{\text{\rm vdW}{\scriptscriptstyle\square}{\mathcal N}}\right)^2$&[5pt]
             \ref{plots:plotHVDW}& $T^2$
             \\};
\end{tikzpicture}
\caption{\small\slshape Rejection frequencies, for  spherical Student samples (2.1 degrees of freedom; see~5.1.2\vspace{.5mm}) and various sample sizes, of Hotelling's test based on $T^2$ and the van der Waerden tests based on $\left(\tenq{T}^{(n)}_{\text{\rm vdW}{\pm}}\right)^2$,  $\left(\tenq{T}^{(n)}_{\text{\rm vdW}{\scriptscriptstyle\square}}\right)^2$, $\left(\tenq{T}^{(n)}_{\text{\rm vdW}{\pm}{\mathcal N}}\right)^2$, 
and~$\left(\tenq{T}^{(n)}_{\text{\rm vdW}{\scriptscriptstyle\square}{\mathcal N}}\right)^2$, respectively,  as  functions of the shift $\eta\vspace{.8mm}$;   $N=500$ replications.\vspace{-6mm}}
\label{fig: BivCompTsphVDW}
\end{figure}\vspace{-5mm}
 \begin{figure}[h]
\centering
\begin{tikzpicture}
    \begin{groupplot}[group style={
                      group name=myplot,
                      group size= 3 by 1
                      },height=5cm,width=4.5cm,
                      ytick={0,0.2,0.4,0.6,0.8,1},
                      ymax=1.1]
\nextgroupplot[title={$n=100$}]\addplot [domain=0:.5, color=black] coordinates { 
(0,0.052)
(0.1,0.202)
(0.2,0.588)
(0.3,0.872000000000001)
(0.4,0.984000000000001)
(0.5,0.990000000000001)
}; 
\label{plots:plot1VDW} 
\addplot [domain=0:.5, color=brightgreen] coordinates { 
(0,0.05)
(0.1,0.172)
(0.2,0.572)
(0.3,0.864000000000001)
(0.4,0.978000000000001)
(0.5,0.994000000000001)
}; 
\label{plots:plot2VDW} 
\addplot [domain=0:.5, color=orange] coordinates { 
(0, 0.072)
(0.1, .1)
(0.2, .25)
(0.3, .508)
(0.4, .758)
(0.5, .886)
}; 
\label{plots:plotHVDW} 
\addplot [domain=0:.5, color=blue] coordinates { 
(0,0.056)
(0.1,0.196)
(0.2,0.592)
(0.3,0.880000000000001)
(0.4,0.982000000000001)
(0.5,0.994000000000001)
}; 
\label{plots:plot3VDW} 
\addplot [domain=0:.5, color=red] coordinates { 
(0,0.05)
(0.1,0.186)
(0.2,0.578)
(0.3,0.860000000000001)
(0.4,0.982000000000001)
(0.5,0.990000000000001)
}; 
\label{plots:plot4VDW} 
\addplot [domain=0:.5, color=grey] coordinates { 
(0,0.05)
(0.1,0.05)
(0.2,0.05)
(0.3,0.05)
(0.4,0.05)
(0.5,0.05)
}; 
\label{plots:plot5VDW} 
\nextgroupplot[title={$n=200$}]\addplot [domain=0:.5, color=black] coordinates { 
(0,0.04)
(0.1,0.35)
(0.2,0.900000000000001)
(0.3,0.994000000000001)
(0.4,1)
(0.5,1)
}; 
\addplot [domain=0:.5, color=brightgreen] coordinates { 
(0,0.042)
(0.1,0.338)
(0.2,0.888000000000001)
(0.3,0.992000000000001)
(0.4,1)
(0.5,1)
}; 
\addplot [domain=0:.5, color=blue] coordinates { 
(0,0.034)
(0.1,0.328)
(0.2,0.886000000000001)
(0.3,0.992000000000001)
(0.4,1)
(0.5,1)
}; 
\addplot [domain=0:.5, color=orange] coordinates { 
(0, 0.04)
(0.1, .156)
(0.2, .444)
(0.3, .816)
(0.4, .976)
(0.5, .998)
}; 
\addplot [domain=0:.5, color=red] coordinates { 
(0,0.036)
(0.1,0.356)
(0.2,0.870000000000001)
(0.3,0.992000000000001)
(0.4,1)
(0.5,1)
}; 
\addplot [domain=0:.5, color=grey] coordinates { 
(0,0.05)
(0.1,0.05)
(0.2,0.05)
(0.3,0.05)
(0.4,0.05)
(0.5,0.05)
}; 
\nextgroupplot[title={$n=400$}]\addplot [domain=0:.5, color=black] coordinates { 
(0,0.058)
(0.1,0.684)
(0.2,0.998000000000001)
(0.3,1)
(0.4,1)
(0.5,1)
}; 
\addplot [domain=0:.5, color=brightgreen] coordinates { 
(0,0.06)
(0.1,0.666)
(0.2,0.998000000000001)
(0.3,1)
(0.4,1)
(0.5,1)
}; 
\addplot [domain=0:.5, color=orange] coordinates { 
(0, 0.046)
(0.1, .248)
(0.2, .768)
(0.3, .984)
(0.4, 1)
(0.5, 1)
}; 
\addplot [domain=0:.5, color=blue] coordinates { 
(0,0.062)
(0.1,0.686)
(0.2,0.998000000000001)
(0.3,1)
(0.4,1)
(0.5,1)
}; 
\addplot [domain=0:.5, color=red] coordinates { 
(0,0.056)
(0.1,0.67)
(0.2,0.996000000000001)
(0.3,1)
(0.4,1)
(0.5,1)
}; 
\addplot [domain=0:.5, color=grey] coordinates { 
(0,0.05)
(0.1,0.05)
(0.2,0.05)
(0.3,0.05)
(0.4,0.05)
(0.5,0.05)
}; 
            \end{groupplot}      
\path (myplot c1r1.south west|-current bounding box.south)--
      coordinate(legendpos)
      (myplot c3r1.south east|-current bounding box.south);
\matrix[
    matrix of nodes,
    anchor=south,
    draw,
    inner sep=0.2em,
    draw
  ]at([yshift=-6ex]legendpos)
  {
     \ref{plots:plot3VDW}& $\left(\tenq{T}^{(n)}_{\text{\rm vdW}{\pm}}\right)^2$&[5pt]
        \ref{plots:plot4VDW}&$\left(\tenq{T}^{(n)}_{\text{\rm vdW}{\scriptscriptstyle\square}}\right)^2$ &[5pt]
    \ref{plots:plot2VDW}& $\left(\tenq{T}^{(n)}_{\text{\rm vdW}{\pm}{\mathcal N}}\right)^2$&[5pt]
             \ref{plots:plot1VDW}& $\left(\tenq{T}^{(n)}_{\text{\rm vdW}{\scriptscriptstyle\square}{\mathcal N}}\right)^2$ &[5pt]
                 \ref{plots:plotHVDW}& $T^2$
            \\};
\end{tikzpicture}
\caption{\small\slshape Rejection frequencies, for  ``banana-shaped'' samples (see~Section~\ref{sec513}) and various sample sizes, of Hotelling's test based on $T^2$ and the van der Waerden tests based on $\left(\tenq{T}^{(n)}_{\text{\rm vdW}{\pm}}\right)^2$,  $\left(\tenq{T}^{(n)}_{\text{\rm vdW}{\scriptscriptstyle\square}}\right)^2\vspace{1mm}$, $\left(\tenq{T}^{(n)}_{\text{\rm vdW}{\pm}{\mathcal N}}\right)^2$, 
and~$\left(\tenq{T}^{(n)}_{\text{\rm vdW}{\scriptscriptstyle\square}{\mathcal N}}\right)^2$, respectively,  as  functions of the shift $\eta\vspace{.6mm}$;   $N=500$ replications.\vspace{-4mm}} 
\label{fig: BivCompBananaVDW}
\end{figure}


\subsubsection{``Banana-shaped'' samples}\label{sec513}\vspace{-3mm}

Same ``banana-shaped'' mixtures as in Section~\ref{sec413}.  
Rejection frequencies 
 over~$500$ 
 replications are
 shown (as  functions of $\eta$) in  Figure~\ref{fig: BivCompBananaVDW}.  The conclusions are quite
 similar as in the previous case: the empirical power curves of the various van der Waerden tests are essentially indistinguishable, while significantly outperforming the Hotelling ones.\vspace{-4mm}


\subsubsection{Samples   with independent Cauchy marginals}\label{sec514}\vspace{-3mm}

Same  samples with independent Cauchy marginals as in Section~\ref{sec414}. 
Rejection frequencies over~$N=500$ replications are
 shown (as  functions of $\eta$) in  Figure~\ref{fig: BivCompCindVDW}.
  Again, all powers are much less than in the Gaussian case, but Hotelling is totally inefficient. \vspace{-4mm}

\begin{figure}[h]
\centering
\begin{tikzpicture}
    \begin{groupplot}[group style={
                      group name=myplot,
                      group size= 3 by 1
                      },height=5cm,width=4.5cm,
                      ytick={0,0.2,0.4,0.6,0.8,1},
                      ymax=1.1]
\nextgroupplot[title={$n=100$}]\addplot [domain=0:.5, color=black] coordinates { 
(0,0.032)
(0.1,0.062)
(0.2,0.068)
(0.3,0.126)
(0.4,0.202)
(0.5,0.272)
}; 
\label{plots:plot1VDW} 
\addplot [domain=0:.5, color=brightgreen] coordinates { 
(0,0.05)
(0.1,0.056)
(0.2,0.07)
(0.3,0.124)
(0.4,0.206)
(0.5,0.286)
}; 
\label{plots:plot2VDW} 
\addplot [domain=0:.5, color=blue] coordinates { 
(0,0.054)
(0.1,0.06)
(0.2,0.0800000000000001)
(0.3,0.13)
(0.4,0.23)
(0.5,0.294)
}; 
\addplot [domain=0:.5, color=orange] coordinates { 
(0, 0.014)
(0.1, .012)
(0.2, .016)
(0.3, .026)
(0.4, .034)
(0.5, .028)
}; 
\label{plots:plotHVDW} 
\label{plots:plot3VDW} 
\addplot [domain=0:.5, color=red] coordinates { 
(0,0.038)
(0.1,0.058)
(0.2,0.064)
(0.3,0.14)
(0.4,0.226)
(0.5,0.278)
}; 
\label{plots:plot4VDW} 
\addplot [domain=0:.5, color=grey] coordinates { 
(0,0.05)
(0.1,0.05)
(0.2,0.05)
(0.3,0.05)
(0.4,0.05)
(0.5,0.05)
}; 
\label{plots:plot5VDW} 
\nextgroupplot[title={$n=200$}]\addplot [domain=0:.5, color=black] coordinates { 
(0,0.06)
(0.1,0.05)
(0.2,0.122)
(0.3,0.236)
(0.4,0.338)
(0.5,0.482)
}; 
\addplot [domain=0:.5, color=brightgreen] coordinates { 
(0,0.058)
(0.1,0.046)
(0.2,0.134)
(0.3,0.236)
(0.4,0.37)
(0.5,0.518)
}; 
\addplot [domain=0:.5, color=blue] coordinates { 
(0,0.054)
(0.1,0.05)
(0.2,0.132)
(0.3,0.234)
(0.4,0.364)
(0.5,0.53)
}; 
\addplot [domain=0:.5, color=orange] coordinates { 
(0, 0.018)
(0.1, .024)
(0.2, .022)
(0.3, .036)
(0.4, .046)
(0.5, .03)
}; 
\addplot [domain=0:.5, color=red] coordinates { 
(0,0.054)
(0.1,0.056)
(0.2,0.136)
(0.3,0.234)
(0.4,0.342)
(0.5,0.506)
}; 
\addplot [domain=0:.5, color=grey] coordinates { 
(0,0.05)
(0.1,0.05)
(0.2,0.05)
(0.3,0.05)
(0.4,0.05)
(0.5,0.05)
}; 
\nextgroupplot[title={$n=400$}]\addplot [domain=0:.5, color=black] coordinates { 
(0,0.038)
(0.1,0.1)
(0.2,0.218)
(0.3,0.402)
(0.4,0.638)
(0.5,0.790000000000001)
}; 
\addplot [domain=0:.5, color=brightgreen] coordinates { 
(0,0.04)
(0.1,0.0820000000000001)
(0.2,0.236)
(0.3,0.43)
(0.4,0.68)
(0.5,0.822000000000001)
}; 
\addplot [domain=0:.5, color=blue] coordinates { 
(0,0.046)
(0.1,0.0940000000000001)
(0.2,0.232)
(0.3,0.418)
(0.4,0.690000000000001)
(0.5,0.832000000000001)
}; 
\addplot [domain=0:.5, color=orange] coordinates { 
(0, 0.028)
(0.1, .026)
(0.2, .012)
(0.3, .014)
(0.4, .02)
(0.5, .038)
}; 
\addplot [domain=0:.5, color=red] coordinates { 
(0,0.038)
(0.1,0.0960000000000001)
(0.2,0.24)
(0.3,0.41)
(0.4,0.66)
(0.5,0.806000000000001)
}; 
\addplot [domain=0:.5, color=grey] coordinates { 
(0,0.05)
(0.1,0.05)
(0.2,0.05)
(0.3,0.05)
(0.4,0.05)
(0.5,0.05)
}; 
            \end{groupplot}      
\path (myplot c1r1.south west|-current bounding box.south)--
      coordinate(legendpos)
      (myplot c3r1.south east|-current bounding box.south);
\matrix[
    matrix of nodes,
    anchor=south,
    draw,
    inner sep=0.2em,
    draw
  ]at([yshift=-6ex]legendpos)
  {
   \ref{plots:plot3VDW}& $\left(\tenq{T}^{(n)}_{\text{\rm vdW}{\pm}}\right)^2$&[5pt]
        \ref{plots:plot4VDW}&$\left(\tenq{T}^{(n)}_{\text{\rm vdW}{\scriptscriptstyle\square}}\right)^2$ &[5pt]
    \ref{plots:plot2VDW}& $\left(\tenq{T}^{(n)}_{\text{\rm vdW}{\pm}{\mathcal N}}\right)^2$&[5pt]
             \ref{plots:plot1VDW}& $\left(\tenq{T}^{(n)}_{\text{\rm vdW}{\scriptscriptstyle\square}{\mathcal N}}\right)^2$ &[5pt]
             \ref{plots:plotHVDW}& $T^2$
              \\};
\end{tikzpicture}
\caption{\small\slshape Rejection frequencies, for samples with independent Cauchy marginals  (see~Section~\ref{sec514}) and various sample sizes, of Hotelling's test based on $T^2$ and the van der Waerden tests based on $\left(\tenq{T}^{(n)}_{\text{\rm vdW}{\pm}}\right)^2$,  $\left(\tenq{T}^{(n)}_{\text{\rm vdW}{\scriptscriptstyle\square}}\right)^2$, $\left(\tenq{T}^{(n)}_{\text{\rm vdW}{\pm}{\mathcal N}}\right)^2$, 
and~$\left(\tenq{T}^{(n)}_{\text{\rm vdW}{\scriptscriptstyle\square}{\mathcal N}}\right)^2$, respectively,  as  functions of the shift $\eta\vspace{1mm}$;   $N=500$ replications.\vspace{-6mm}}  
\label{fig: BivCompCindVDW}
\end{figure}

\subsubsection{Nonspherical Gaussian samples}\label{sec515}\vspace{-3mm}

Same  correlated Gaussian samples  as in Section~\ref{sec415}. 
Rejection frequencies over~$N=~\!500$ replications are
 shown (as  functions of $\eta$) in  Figure~\ref{fig: BivCompGausCorVDW}. 
  The non-specification of the covariance matrix apparently has no impact on the comparative performance of Hotelling and its rank-based van der Waerden competitors, which all coincide. \vspace{-4mm}

\begin{figure}[h]
\centering
\begin{tikzpicture}
    \begin{groupplot}[group style={
                      group name=myplot,
                      group size= 3 by 1
                      },height=5cm,width=4.5cm,
                      ytick={0,0.2,0.4,0.6,0.8,1},
                      ymax=1.1]
\nextgroupplot[title={$n=100$}]\addplot [domain=0:.5, color=black] coordinates { 
(0,0.05)
(0.1,0.0740000000000001)
(0.2,0.184)
(0.3,0.47)
(0.4,0.624)
(0.5,0.846000000000001)
}; 
\label{plots:plot1VDW} 
\addplot [domain=0:.5, color=brightgreen] coordinates { 
(0,0.04)
(0.1,0.0720000000000001)
(0.2,0.196)
(0.3,0.466)
(0.4,0.644)
(0.5,0.850000000000001)
}; 
\label{plots:plot2VDW} 
\addplot [domain=0:.5, color=blue] coordinates { 
(0,0.046)
(0.1,0.0820000000000001)
(0.2,0.212)
(0.3,0.482)
(0.4,0.638)
(0.5,0.854000000000001)
}; 
\label{plots:plot3VDW} 
\addplot [domain=0:.5, color=orange] coordinates { 
(0,0.07)
(0.1,0.056)
(0.2,0.122)
(0.3,0.244)
(0.4,0.434)
(0.5,0.628)
}; 
\label{plots:plotHVDW} 
\addplot [domain=0:.5, color=red] coordinates { 
(0,0.042)
(0.1,0.0740000000000001)
(0.2,0.198)
(0.3,0.462)
(0.4,0.63)
(0.5,0.832000000000001)
}; 
\label{plots:plot4VDW} 
\addplot [domain=0:.5, color=grey] coordinates { 
(0,0.05)
(0.1,0.05)
(0.2,0.05)
(0.3,0.05)
(0.4,0.05)
(0.5,0.05)
}; 
\label{plots:plot5VDW} 
\nextgroupplot[title={$n=200$}]\addplot [domain=0:.5, color=black] coordinates { 
(0,0.052)
(0.1,0.104)
(0.2,0.374)
(0.3,0.770000000000001)
(0.4,0.926000000000001)
(0.5,0.992000000000001)
}; 
\addplot [domain=0:.5, color=brightgreen] coordinates { 
(0,0.056)
(0.1,0.114)
(0.2,0.37)
(0.3,0.780000000000001)
(0.4,0.924000000000001)
(0.5,0.994000000000001)
}; 
\addplot [domain=0:.5, color=blue] coordinates { 
(0,0.044)
(0.1,0.102)
(0.2,0.38)
(0.3,0.780000000000001)
(0.4,0.930000000000001)
(0.5,0.994000000000001)
}; 
\addplot [domain=0:.5, color=orange] coordinates { 
(0,0.05)
(0.1,0.104)
(0.2,0.30)
(0.3,0.486)
(0.4,0.73)
(0.5,0.902)
}; 
\addplot [domain=0:.5, color=red] coordinates { 
(0,0.054)
(0.1,0.114)
(0.2,0.37)
(0.3,0.782000000000001)
(0.4,0.924000000000001)
(0.5,0.994000000000001)
}; 
\addplot [domain=0:.5, color=grey] coordinates { 
(0,0.05)
(0.1,0.05)
(0.2,0.05)
(0.3,0.05)
(0.4,0.05)
(0.5,0.05)
}; 
\nextgroupplot[title={$n=400$}]\addplot [domain=0:.5, color=black] coordinates { 
(0,0.042)
(0.1,0.218)
(0.2,0.690000000000001)
(0.3,0.968000000000001)
(0.4,1)
(0.5,1)
}; 
\addplot [domain=0:.5, color=brightgreen] coordinates { 
(0,0.044)
(0.1,0.214)
(0.2,0.704000000000001)
(0.3,0.966000000000001)
(0.4,1)
(0.5,1)
}; 
\addplot [domain=0:.5, color=blue] coordinates { 
(0,0.044)
(0.1,0.218)
(0.2,0.712000000000001)
(0.3,0.968000000000001)
(0.4,1)
(0.5,1)
}; 
\addplot [domain=0:.5, color=orange] coordinates { 
(0,0.064)
(0.1,0.146)
(0.2,0.482)
(0.3,0.826)
(0.4,0.978)
(0.5,0.998)
}; 
\addplot [domain=0:.5, color=red] coordinates { 
(0,0.044)
(0.1,0.218)
(0.2,0.704000000000001)
(0.3,0.972000000000001)
(0.4,1)
(0.5,1)
}; 
\addplot [domain=0:.5, color=grey] coordinates { 
(0,0.05)
(0.1,0.05)
(0.2,0.05)
(0.3,0.05)
(0.4,0.05)
(0.5,0.05)
}; 
            \end{groupplot}      
\path (myplot c1r1.south west|-current bounding box.south)--
      coordinate(legendpos)
      (myplot c3r1.south east|-current bounding box.south);
\matrix[
    matrix of nodes,
    anchor=south,
    draw,
    inner sep=0.2em,
    draw
  ]at([yshift=-6ex]legendpos)
  {
      \ref{plots:plot3VDW}& $\left(\tenq{T}^{(n)}_{\text{\rm vdW}{\pm}}\right)^2$&[5pt]
        \ref{plots:plot4VDW}&$\left(\tenq{T}^{(n)}_{\text{\rm vdW}{\scriptscriptstyle\square}}\right)^2$ &[5pt]
    \ref{plots:plot2VDW}& $\left(\tenq{T}^{(n)}_{\text{\rm vdW}{\pm}{\mathcal N}}\right)^2$&[5pt]
             \ref{plots:plot1VDW}& $\left(\tenq{T}^{(n)}_{\text{\rm vdW}{\scriptscriptstyle\square}{\mathcal N}}\right)^2$ &[5pt]
              \ref{plots:plotHVDW}& $T^2$
 \\};
\end{tikzpicture}
\caption{\small\slshape Rejection frequencies, for  nonspherical Gaussian samples (see~Section~\ref{sec515}) and various sample sizes, of Hotelling's test based on $T^2$ and the van der Waerden tests based on $\left(\tenq{T}^{(n)}_{\text{\rm vdW}{\pm}}\right)^2\vspace{-.2mm}$,  $\left(\tenq{T}^{(n)}_{\text{\rm vdW}{\scriptscriptstyle\square}}\right)^2$, $\left(\tenq{T}^{(n)}_{\text{\rm vdW}{\pm}{\mathcal N}}\right)^2$, 
and~$\left(\tenq{T}^{(n)}_{\text{\rm vdW}{\scriptscriptstyle\square}{\mathcal N}}\right)^2$, respectively,  as  functions of the shift $\eta\vspace{1mm}$;   $N=500$ replications.\vspace{-1mm}} 
\label{fig: BivCompGausCorVDW}
\end{figure}

\subsubsection{Spherical Cauchy samples}\label{sec516}\vspace{-3mm}

Same  spherical Cauchy samples   as in Section~\ref{sec416}. 
Rejection frequencies over~$500$ replications are
 shown (as  functions of $\eta$) in  Figure~\ref{fig: BivCompCsphVDW}. All tests perform similarly except, of course, for Hotelling, which fails completely. \vspace{-4mm}

\begin{figure}[h]
\centering
\begin{tikzpicture}
    \begin{groupplot}[group style={
                      group name=myplot,
                      group size= 3 by 1
                      },height=5cm,width=4.5cm,
                      ytick={0,0.2,0.4,0.6,0.8,1},
                      ymax=1.1]
\nextgroupplot[title={$n=100$}]\addplot [domain=0:.5, color=black] coordinates { 
(0,0.05)
(0.1,0.052)
(0.2,0.0880000000000001)
(0.3,0.152)
(0.4,0.28)
(0.5,0.342)
}; 
\label{plots:plot1VDW} 
\addplot [domain=0:.5, color=brightgreen] coordinates { 
(0,0.048)
(0.1,0.056)
(0.2,0.0860000000000001)
(0.3,0.156)
(0.4,0.272)
(0.5,0.362)
}; 
\label{plots:plot2VDW} 
\addplot [domain=0:.5, color=blue] coordinates { 
(0,0.044)
(0.1,0.05)
(0.2,0.0960000000000001)
(0.3,0.166)
(0.4,0.286)
(0.5,0.374)
}; 
\label{plots:plot3VDW} 
\addplot [domain=0:.5, color=orange] coordinates { 
(0,0.01)
(0.1,0.018)
(0.2,0.022)
(0.3,0.028)
(0.4,0.048)
(0.5,0.04)
}; 
\label{plots:plotH} 
\addplot [domain=0:.5, color=red] coordinates { 
(0,0.048)
(0.1,0.064)
(0.2,0.0860000000000001)
(0.3,0.156)
(0.4,0.266)
(0.5,0.354)
}; 
\label{plots:plot4VDW} 
\addplot [domain=0:.5, color=grey] coordinates { 
(0,0.05)
(0.1,0.05)
(0.2,0.05)
(0.3,0.05)
(0.4,0.05)
(0.5,0.05)
}; 
\label{plots:plot5VDW} 
\nextgroupplot[title={$n=200$}]\addplot [domain=0:.5, color=black] coordinates { 
(0,0.048)
(0.1,0.0740000000000001)
(0.2,0.158)
(0.3,0.266)
(0.4,0.436)
(0.5,0.628)
}; 
\addplot [domain=0:.5, color=brightgreen] coordinates { 
(0,0.042)
(0.1,0.0840000000000001)
(0.2,0.182)
(0.3,0.308)
(0.4,0.45)
(0.5,0.654)
}; 
\addplot [domain=0:.5, color=blue] coordinates { 
(0,0.038)
(0.1,0.0880000000000001)
(0.2,0.162)
(0.3,0.29)
(0.4,0.458)
(0.5,0.656)
}; 
\addplot [domain=0:.5, color=orange] coordinates { 
(0,0.012)
(0.1,0.018)
(0.2,0.016)
(0.3,0.022)
(0.4,0.03)
(0.5,0.048)
}; 
\addplot [domain=0:.5, color=red] coordinates { 
(0,0.05)
(0.1,0.0760000000000001)
(0.2,0.156)
(0.3,0.276)
(0.4,0.44)
(0.5,0.638)
}; 
\addplot [domain=0:.5, color=grey] coordinates { 
(0,0.05)
(0.1,0.05)
(0.2,0.05)
(0.3,0.05)
(0.4,0.05)
(0.5,0.05)
}; 
\nextgroupplot[title={$n=400$}]\addplot [domain=0:.5, color=black] coordinates { 
(0,0.0740000000000001)
(0.1,0.0800000000000001)
(0.2,0.21)
(0.3,0.486)
(0.4,0.738000000000001)
(0.5,0.914000000000001)
}; 
\addplot [domain=0:.5, color=brightgreen] coordinates { 
(0,0.0720000000000001)
(0.1,0.0840000000000001)
(0.2,0.218)
(0.3,0.518)
(0.4,0.760000000000001)
(0.5,0.934000000000001)
}; 
\addplot [domain=0:.5, color=blue] coordinates { 
(0,0.0760000000000001)
(0.1,0.0840000000000001)
(0.2,0.216)
(0.3,0.522)
(0.4,0.772000000000001)
(0.5,0.942000000000001)
}; 
\addplot [domain=0:.5, color=orange] coordinates { 
(0,0.032)
(0.1,0.014)
(0.2,0.022)
(0.3,0.014)
(0.4,0.02400000001)
(0.5,0.032000000000001)
}; 
\addplot [domain=0:.5, color=red] coordinates { 
(0,0.0760000000000001)
(0.1,0.0920000000000001)
(0.2,0.226)
(0.3,0.508)
(0.4,0.736000000000001)
(0.5,0.920000000000001)
}; 
\addplot [domain=0:.5, color=grey] coordinates { 
(0,0.05)
(0.1,0.05)
(0.2,0.05)
(0.3,0.05)
(0.4,0.05)
(0.5,0.05)
}; 

            \end{groupplot}      
\path (myplot c1r1.south west|-current bounding box.south)--
      coordinate(legendpos)
      (myplot c3r1.south east|-current bounding box.south);
\matrix[
    matrix of nodes,
    anchor=south,
    draw,
    inner sep=0.2em,
    draw
  ]at([yshift=-6ex]legendpos)
  {
     \ref{plots:plot3VDW}& $\left(\tenq{T}^{(n)}_{\text{\rm vdW}{\pm}}\right)^2$&[5pt]
        \ref{plots:plot4VDW}&$\left(\tenq{T}^{(n)}_{\text{\rm vdW}{\scriptscriptstyle\square}}\right)^2$ &[5pt]
    \ref{plots:plot2VDW}& $\left(\tenq{T}^{(n)}_{\text{\rm vdW}{\pm}{\mathcal N}}\right)^2$&[5pt]
             \ref{plots:plot1VDW}& $\left(\tenq{T}^{(n)}_{\text{\rm vdW}{\scriptscriptstyle\square}{\mathcal N}}\right)^2$ &[5pt]
                \ref{plots:plotHVDW}& $T^2$
 \\};
\end{tikzpicture}
\caption{\small\slshape Rejection frequencies,  for  spherical Cauchy samples (see~Section~\ref{sec516}) and various sample sizes, of Hotelling's test based on $T^2$ and the van der Waerden tests based on $\left(\tenq{T}^{(n)}_{\text{\rm vdW}{\pm}}\right)^2$,  $\left(\tenq{T}^{(n)}_{\text{\rm vdW}{\scriptscriptstyle\square}}\right)^2\vspace{0mm}$, $\left(\tenq{T}^{(n)}_{\text{\rm vdW}{\pm}{\mathcal N}}\right)^2$, 
and~$\left(\tenq{T}^{(n)}_{\text{\rm vdW}{\scriptscriptstyle\square}{\mathcal N}}\right)^2$, respectively,  as  functions of the shift $\eta$;   $N=500$ replications.\vspace{-4mm}} 
\label{fig: BivCompCsphVDW}
\end{figure}

\subsection{van der Waerden-type statistics in dimension $d=5$}



\subsubsection{Spherical Gaussian samples}\label{sec521}\vspace{-3mm}

Same  spherical Gaussian samples   as in Section~\ref{sec421}. 
Rejection frequencies over~$N=~\!500$ replications are
 shown (as  functions of $\eta$) in  Figure~\ref{fig: CompGausVDW5}. A ``small sample'' superiority of  Hotelling for $n=200$ rapidly disappears as $n$ increases; all van der Waerden tests yield the same performance. \vspace{-4mm}

\begin{figure}[h]
\centering
\begin{tikzpicture}
    \begin{groupplot}[group style={
                      group name=myplot,
                      group size= 3 by 1
                      },height=5cm,width=4.5cm,
                      ytick={0,0.2,0.4,0.6,0.8,1},
                      ymax=1.1]
\nextgroupplot[title={$n=200$}]\addplot [domain=0:.5, color=black] coordinates { 
(0,0.048)
(0.1,0.152)
(0.2,0.474)
(0.3,0.878000000000001)
(0.4,0.990000000000001)
(0.5,1)
}; 
\addplot [domain=0:.5, color=brightgreen] coordinates { 
(0,0.044)
(0.1,0.158)
(0.2,0.562)
(0.3,0.922000000000001)
(0.4,0.996000000000001)
(0.5,1)
}; 
\addplot [domain=0:.5, color=blue] coordinates { 
(0,0.056)
(0.1,0.17)
(0.2,0.512)
(0.3,0.904000000000001)
(0.4,0.994000000000001)
(0.5,1)
}; 
\addplot [domain=0:.5, color=orange] coordinates { 
(0,0.064)
(0.1,0.188)
(0.2,0.676)
(0.3,0.966000000000001)
(0.4,1)
(0.5,1)
}; 
\addplot [domain=0:.5, color=red] coordinates { 
(0,0.046)
(0.1,0.166)
(0.2,0.474)
(0.3,0.894000000000001)
(0.4,0.994000000000001)
(0.5,1)
}; 
\addplot [domain=0:.5, color=grey] coordinates { 
(0,0.05)
(0.1,0.05)
(0.2,0.05)
(0.3,0.05)
(0.4,0.05)
(0.5,0.05)
}; 
\nextgroupplot[title={$n=400$}]\addplot [domain=0:.5, color=black] coordinates { 
(0,0.068)
(0.1,0.32)
(0.2,0.872000000000001)
(0.3,1)
(0.4,1)
(0.5,1)
}; 
\addplot [domain=0:.5, color=brightgreen] coordinates { 
(0,0.0720000000000001)
(0.1,0.294)
(0.2,0.872000000000001)
(0.3,0.998000000000001)
(0.4,1)
(0.5,1)
}; 
\addplot [domain=0:.5, color=blue] coordinates { 
(0,0.0760000000000001)
(0.1,0.284)
(0.2,0.848000000000001)
(0.3,1)
(0.4,1)
(0.5,1)
}; 
\addplot [domain=0:.5, color=orange] coordinates { 
(0,0.056)
(0.1,0.4)
(0.2,0.952)
(0.3,1)
(0.4,1)
(0.5,1)
}; 
\addplot [domain=0:.5, color=red] coordinates { 
(0,0.068)
(0.1,0.284)
(0.2,0.876000000000001)
(0.3,1)
(0.4,1)
(0.5,1)
}; 
\addplot [domain=0:.5, color=grey] coordinates { 
(0,0.05)
(0.1,0.05)
(0.2,0.05)
(0.3,0.05)
(0.4,0.05)
(0.5,0.05)
}; 
\nextgroupplot[title={$n=800$}]\addplot [domain=0:.5, color=black] coordinates { 
(0,0.044)
(0.05,0.154)
(0.1,0.574)
(0.15,0.950000000000001)
(0.2,0.998000000000001)
(0.25,1)
(0.5,1)
}; 
\label{plots:plot1VDW5D5} 
\addplot [domain=0:.5, color=brightgreen] coordinates { 
(0,0.03)
(0.05,0.156)
(0.1,0.564)
(0.15,0.932000000000001)
(0.2,0.994000000000001)
(0.25,1)
(0.5,1)
}; 
\label{plots:plot2VDW5D5} 
\addplot [domain=0:.5, color=blue] coordinates { 
(0,0.04)
(0.05,0.156)
(0.1,0.558)
(0.15,0.934000000000001)
(0.2,0.996000000000001)
(0.25,1)
(0.5,1)
}; 
\label{plots:plot3VDW5D5} 
\addplot [domain=0:.5, color=red] coordinates { 
(0,0.044)
(0.05,0.16)
(0.1,0.582)
(0.15,0.924000000000001)
(0.2,0.996000000000001)
(0.25,1)
(0.5,1)
}; 
\label{plots:plot4VDW5D5} 
\addplot [domain=0:.5, color=orange] coordinates { 
(0,0.038)
(0.05,0.194)
(0.1,0.68)
(0.15,0.968000000000001)
(0.2,0.998000000000001)
(0.25,1)
(0.5,1)
}; 
\label{plots:plot5VDW5D5} 
\addplot [domain=0:.5, color=grey] coordinates { 
(0,0.05)
(0.05,0.05)
(0.1,0.05)
(0.15,0.05)
(0.2,0.05)
(0.25,0.05)
(0.5,0.05)
}; 
\label{plots:plot6VDW5D5} 

            \end{groupplot}      
\path (myplot c1r1.south west|-current bounding box.south)--
      coordinate(legendpos)
      (myplot c3r1.south east|-current bounding box.south);
\matrix[
    matrix of nodes,
    anchor=south,
    draw,
    inner sep=0.2em,
    draw
  ]at([yshift=-6ex]legendpos)
  {
     \ref{plots:plot3VDW}& $\left(\tenq{T}^{(n)}_{\text{\rm vdW}{\pm}}\right)^2$&[5pt]
        \ref{plots:plot4VDW}&$\left(\tenq{T}^{(n)}_{\text{\rm vdW}{\scriptscriptstyle\square}}\right)^2$ &[5pt]
    \ref{plots:plot2VDW}& $\left(\tenq{T}^{(n)}_{\text{\rm vdW}{\pm}{\mathcal N}}\right)^2$&[5pt]
             \ref{plots:plot1VDW}& $\left(\tenq{T}^{(n)}_{\text{\rm vdW}{\scriptscriptstyle\square}{\mathcal N}}\right)^2$&[5pt]
 \ref{plots:plotHVDW}& $T^2$
 \\};
\end{tikzpicture}
\caption{\small\slshape Rejection frequencies, for samples with 5-dimensional spherical Gaussian distributions (see~Section~\ref{sec521}) and various sample sizes, of Hotelling's test based on $T^2$ and the van der Waerden tests based on $\left(\tenq{T}^{(n)}_{\text{\rm vdW}{\pm}}\right)^2$,  $\left(\tenq{T}^{(n)}_{\text{\rm vdW}{\scriptscriptstyle\square}}\right)^2\vspace{1mm}$, $\left(\tenq{T}^{(n)}_{\text{\rm vdW}{\pm}{\mathcal N}}\right)^2$, 
and~$\left(\tenq{T}^{(n)}_{\text{\rm vdW}{\scriptscriptstyle\square}{\mathcal N}}\right)^2\vspace{.4mm}$, respectively,  as  functions of the shift $\eta$;   $N=500$ replications.} 
\label{fig: CompGausVDW5}\vspace{0mm}
\end{figure}



\subsubsection{Spherical Student samples}\label{sec522}\vspace{-3mm}

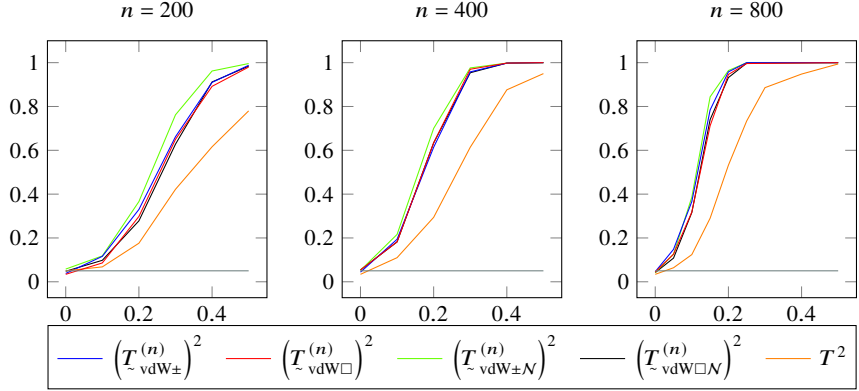
\begin{figure}[h]
\centering
\begin{tikzpicture}
    \begin{groupplot}[group style={
                      group name=myplot,
                      group size= 3 by 1
                      },height=5cm,width=4.5cm,
                      ytick={0,0.2,0.4,0.6,0.8,1},
                      ymax=1.1]
\nextgroupplot[title={$n=200$}]\addplot [domain=0:.5, color=black] coordinates { 
(0,0.048)
(0.1,0.0980000000000001)
(0.2,0.278)
(0.3,0.628)
(0.4,0.912000000000001)
(0.5,0.984000000000001)
}; 
\addplot [domain=0:.5, color=brightgreen] coordinates { 
(0,0.058)
(0.1,0.118)
(0.2,0.366)
(0.3,0.762000000000001)
(0.4,0.962000000000001)
(0.5,0.996000000000001)
}; 
\addplot [domain=0:.5, color=blue] coordinates { 
(0,0.038)
(0.1,0.116)
(0.2,0.33)
(0.3,0.662)
(0.4,0.910000000000001)
(0.5,0.988000000000001)
}; 
\addplot [domain=0:.5, color=orange] coordinates { 
(0,0.05)
(0.1,0.068)
(0.2,0.176)
(0.3,0.422)
(0.4,0.616)
(0.5,0.78)
}; 
\addplot [domain=0:.5, color=red] coordinates { 
(0,0.034)
(0.1,0.0860000000000001)
(0.2,0.298)
(0.3,0.65)
(0.4,0.892000000000001)
(0.5,0.980000000000001)
}; 
\addplot [domain=0:.5, color=grey] coordinates { 
(0,0.05)
(0.1,0.05)
(0.2,0.05)
(0.3,0.05)
(0.4,0.05)
(0.5,0.05)
}; 
\nextgroupplot[title={$n=400$}]\addplot [domain=0:.5, color=black] coordinates { 
(0,0.054)
(0.1,0.182)
(0.2,0.634)
(0.3,0.954000000000001)
(0.4,0.998000000000001)
(0.5,1)
}; 
\addplot [domain=0:.5, color=brightgreen] coordinates { 
(0,0.052)
(0.1,0.216)
(0.2,0.700000000000001)
(0.3,0.976000000000001)
(0.4,0.998000000000001)
(0.5,1)
}; 
\addplot [domain=0:.5, color=blue] coordinates { 
(0,0.044)
(0.1,0.194)
(0.2,0.614)
(0.3,0.958000000000001)
(0.4,0.998000000000001)
(0.5,1)
}; 
\addplot [domain=0:.5, color=orange] coordinates { 
(0,0.034)
(0.1,0.11)
(0.2,0.294)
(0.3,0.614)
(0.4,0.876)
(0.5,0.95)
}; 
\addplot [domain=0:.5, color=red] coordinates { 
(0,0.054)
(0.1,0.182)
(0.2,0.634)
(0.3,0.970000000000001)
(0.4,0.998000000000001)
(0.5,1)
}; 
\addplot [domain=0:.5, color=grey] coordinates { 
(0,0.05)
(0.1,0.05)
(0.2,0.05)
(0.3,0.05)
(0.4,0.05)
(0.5,0.05)
}; 
\nextgroupplot[title={$n=800$}]\addplot [domain=0:.5, color=black] coordinates { 
(0,0.042)
(0.05,0.108)
(0.1,0.316)
(0.15,0.74)
(0.2,0.932000000000001)
(0.25,0.998000000000001)
(0.5,1)
}; 
\label{plots:plot1VDW5D6} 
\addplot [domain=0:.5, color=brightgreen] coordinates { 
(0,0.042)
(0.05,0.126)
(0.1,0.386)
(0.15,0.844000000000001)
(0.2,0.964000000000001)
(0.25,1)
(0.5,1)
}; 
\label{plots:plot2VDW5D6} 
\addplot [domain=0:.5, color=blue] coordinates { 
(0,0.046)
(0.05,0.148)
(0.1,0.366)
(0.15,0.784000000000001)
(0.2,0.958000000000001)
(0.25,1)
(0.5,1)
}; 
\label{plots:plot3VDW5D6} 
\addplot [domain=0:.5, color=red] coordinates { 
(0,0.044)
(0.05,0.126)
(0.1,0.318)
(0.15,0.712)
(0.2,0.946000000000001)
(0.25,0.996000000000001)
(0.5,1)
}; 
\label{plots:plot4VDW5D6} 
\addplot [domain=0:.5, color=orange] coordinates { 
(0,0.034)
(0.05,0.064)
(0.1,0.124)
(0.15,0.29)
(0.2,0.532)
(0.25,0.734)
(0.3, 0.886)
(0.4,0.948)
(0.5,0.994)
}; 
\label{plots:plot5VDW5D6} 
\addplot [domain=0:.5, color=grey] coordinates { 
(0,0.05)
(0.05,0.05)
(0.1,0.05)
(0.15,0.05)
(0.2,0.05)
(0.25,0.05)
(0.5,0.05)
}; 
\label{plots:plot6VDW5D6} 
            \end{groupplot}      
\path (myplot c1r1.south west|-current bounding box.south)--
      coordinate(legendpos)
      (myplot c3r1.south east|-current bounding box.south);
\matrix[
    matrix of nodes,
    anchor=south,
    draw,
    inner sep=0.2em,
    draw
  ]at([yshift=-6ex]legendpos)
  {
     \ref{plots:plot3VDW}& $\left(\tenq{T}^{(n)}_{\text{\rm vdW}{\pm}}\right)^2$&[5pt]
        \ref{plots:plot4VDW}&$\left(\tenq{T}^{(n)}_{\text{\rm vdW}{\scriptscriptstyle\square}}\right)^2$ &[5pt]
    \ref{plots:plot2VDW}& $\left(\tenq{T}^{(n)}_{\text{\rm vdW}{\pm}{\mathcal N}}\right)^2$&[5pt]
             \ref{plots:plot1VDW}& $\left(\tenq{T}^{(n)}_{\text{\rm vdW}{\scriptscriptstyle\square}{\mathcal N}}\right)^2$&[5pt]
             \ref{plots:plotHVDW}& $T^2$
 \\};
\end{tikzpicture}
\caption{\small\slshape Rejection frequencies, for samples with 5-dimensional spherical Student  distributions (2.1 degrees of freedom; see~Section~\ref{sec522}) and various sample sizes, of Hotelling's test based on~$T^2$ and the van der Waerden tests based on $\left(\tenq{T}^{(n)}_{\text{\rm vdW}{\pm}}\right)^2$,  $\left(\tenq{T}^{(n)}_{\text{\rm vdW}{\scriptscriptstyle\square}}\right)^2\vspace{1mm}$, $\left(\tenq{T}^{(n)}_{\text{\rm vdW}{\pm}{\mathcal N}}\right)^2$, 
and~$\left(\tenq{T}^{(n)}_{\text{\rm vdW}{\scriptscriptstyle\square}{\mathcal N}}\right)^2$, respectively,  as  functions of the shift $\eta$;   $N=500$ replications.} 
\label{fig: CompTsphVDW5}
\end{figure}

Same  spherical Student samples   as in Section~\ref{sec422}. 
Rejection frequencies over~$500$ replications are
 shown (as  functions of $\eta$) in  Figure~\ref{fig: CompTsphVDW5}. All van der Waerden tests roughly yield the same performance, with a slight adavantage in favor of~$\left(\tenq{T}^{(n)}_{\text{\rm vdW}{\pm}{\mathcal N}}\right)^2$ in ``small samples.'' \vspace{-4mm}

\subsubsection{Nonspherical Gaussian samples}\label{sec523}\vspace{-3mm}

Same  spherical Gaussian samples   as in Section~\ref{sec423}. 
Rejection frequencies over~$N=~\!500$ replications are
 shown (as  functions of $\eta$) in  Figure~\ref{fig: CompGausCorrVDW5}. The slight advantage for $n=200$ of Hotelling under spherical Gaussian has almost disappeared. All van der Waerden tests yield similar performance.   \vspace{-2mm}

\begin{figure}[h]
\centering
\begin{tikzpicture}
    \begin{groupplot}[group style={
                      group name=myplot,
                      group size= 3 by 1
                      },height=5cm,width=4.5cm,
                      ytick={0,0.2,0.4,0.6,0.8,1},
                      ymax=1.1]
\nextgroupplot[title={$n=200$}]\addplot [domain=0:.5, color=black] coordinates { 
(0,0.054)
(0.1,0.0780000000000001)
(0.2,0.174)
(0.3,0.418)
(0.4,0.688000000000001)
(0.5,0.892000000000001)
}; 
\addplot [domain=0:.5, color=brightgreen] coordinates { 
(0,0.056)
(0.1,0.0880000000000001)
(0.2,0.214)
(0.3,0.476)
(0.4,0.738000000000001)
(0.5,0.926000000000001)
}; 
\addplot [domain=0:.5, color=blue] coordinates { 
(0,0.048)
(0.1,0.066)
(0.2,0.196)
(0.3,0.44)
(0.4,0.690000000000001)
(0.5,0.910000000000001)
}; 
\addplot [domain=0:.5, color=orange] coordinates { 
(0,0.042)
(0.1,0.0780000000000001)
(0.2,0.218)
(0.3,0.52)
(0.4,0.81)
(0.5,0.958)
}; 
\addplot [domain=0:.5, color=red] coordinates { 
(0,0.046)
(0.1,0.06)
(0.2,0.176)
(0.3,0.388)
(0.4,0.662)
(0.5,0.880000000000001)
}; 
\addplot [domain=0:.5, color=grey] coordinates { 
(0,0.05)
(0.1,0.05)
(0.2,0.05)
(0.3,0.05)
(0.4,0.05)
(0.5,0.05)
}; 
\nextgroupplot[title={$n=400$}]\addplot [domain=0:.5, color=black] coordinates { 
(0,0.066)
(0.1,0.11)
(0.2,0.394)
(0.3,0.804000000000001)
(0.4,0.968000000000001)
(0.5,1)
}; 
\addplot [domain=0:.5, color=brightgreen] coordinates { 
(0,0.048)
(0.1,0.118)
(0.2,0.4)
(0.3,0.828000000000001)
(0.4,0.972000000000001)
(0.5,1)
}; 
\addplot [domain=0:.5, color=blue] coordinates { 
(0,0.064)
(0.1,0.11)
(0.2,0.378)
(0.3,0.784000000000001)
(0.4,0.972000000000001)
(0.5,1)
}; 
\addplot [domain=0:.5, color=orange] coordinates { 
(0,0.058)
(0.1,0.126)
(0.2,0.488)
(0.3,0.852000000000001)
(0.4,0.996)
(0.5,1)
}; 
\addplot [domain=0:.5, color=red] coordinates { 
(0,0.062)
(0.1,0.138)
(0.2,0.382)
(0.3,0.784000000000001)
(0.4,0.972000000000001)
(0.5,1)
}; 
\addplot [domain=0:.5, color=grey] coordinates { 
(0,0.05)
(0.1,0.05)
(0.2,0.05)
(0.3,0.05)
(0.4,0.05)
(0.5,0.05)
}; 
\nextgroupplot[title={$n=800$}]\addplot [domain=0:.5, color=black] coordinates { 
(0,0.048)
(0.05,0.0820000000000001)
(0.1,0.214)
(0.15,0.434)
(0.2,0.768000000000001)
(0.25,0.944000000000001)
(0.3, 0.99)
(0.4, 1)
(0.5, 1)
}; 
\label{plots:plot1VDW5D7} 
\addplot [domain=0:.5, color=brightgreen] coordinates { 
(0,0.052)
(0.05,0.0800000000000001)
(0.1,0.252)
(0.15,0.46)
(0.2,0.792000000000001)
(0.25,0.934000000000001)
(0.3, 0.988)
(0.4, 1)
(0.5, 1)
}; 
\label{plots:plot2VDW5D7} 
\addplot [domain=0:.5, color=blue] coordinates { 
(0,0.056)
(0.05,0.072)
(0.1,0.26)
(0.15,0.416)
(0.2,0.778000000000001)
(0.25,0.932000000000001)
(0.3, 0.99)
(0.4, 1)
(0.5, 1)
}; 
\label{plots:plot3VDW5D7} 
\addplot [domain=0:.5, color=red] coordinates { 
(0,0.048)
(0.05,0.078)
(0.1,0.21)
(0.15,0.436)
(0.2,0.754)
(0.25,0.938000000000001)
(0.3, 0.988)
(0.4, 1)
(0.5, 1)
}; 
\label{plots:plot4VDW5D7} 
\addplot [domain=0:.5, color=orange] coordinates { 
(0,0.044)
(0.05,0.0840000000000001)
(0.1,0.264)
(0.15,0.498)
(0.2,0.824000000000001)
(0.25,0.954000000000001)
}; 
\label{plots:plot5VDW5D7} 
\addplot [domain=0:.5, color=grey] coordinates { 
(0,0.05)
(0.05,0.05)
(0.1,0.05)
(0.15,0.05)
(0.2,0.05)
(0.25,0.05)
(0.5,0.05)
}; 
\label{plots:plot6VDW5D7}

            \end{groupplot}      
\path (myplot c1r1.south west|-current bounding box.south)--
      coordinate(legendpos)
      (myplot c3r1.south east|-current bounding box.south);
\matrix[
    matrix of nodes,
    anchor=south,
    draw,
    inner sep=0.2em,
    draw
  ]at([yshift=-6ex]legendpos)
  {
     \ref{plots:plot3VDW}& $\left(\tenq{T}^{(n)}_{\text{\rm vdW}{\pm}}\right)^2$&[5pt]
        \ref{plots:plot4VDW}&$\left(\tenq{T}^{(n)}_{\text{\rm vdW}{\scriptscriptstyle\square}}\right)^2$ &[5pt]
    \ref{plots:plot2VDW}& $\left(\tenq{T}^{(n)}_{\text{\rm vdW}{\pm}{\mathcal N}}\right)^2$&[5pt]
             \ref{plots:plot1VDW}& $\left(\tenq{T}^{(n)}_{\text{\rm vdW}{\scriptscriptstyle\square}{\mathcal N}}\right)^2$&[5pt]
             \ref{plots:plotHVDW}& $T^2$
 \\};
\end{tikzpicture}
\caption{\small\slshape Rejection frequencies,   for nonspherical 5-dimensional  Gaussian samples (see~Section~\ref{sec523}) and various sample sizes, of Hotelling's test based on~$T^2$ and the van der Waerden tests based on $\left(\tenq{T}^{(n)}_{\text{\rm vdW}{\pm}}\right)^2$,  $\left(\tenq{T}^{(n)}_{\text{\rm vdW}{\scriptscriptstyle\square}}\right)^2$, $\left(\tenq{T}^{(n)}_{\text{\rm vdW}{\pm}{\mathcal N}}\right)^2$, 
and~$\left(\tenq{T}^{(n)}_{\text{\rm vdW}{\scriptscriptstyle\square}{\mathcal N}}\right)^2\vspace{1mm}$, respectively,  as  functions of the shift $\eta$;   $N=500$ replications.\vspace{-5mm}} 
\label{fig: CompGausCorrVDW5}
\end{figure}

\begin{figure}[h]
\centering
\begin{tikzpicture}
    \begin{groupplot}[group style={
                      group name=myplot,
                      group size= 3 by 1
                      },height=5cm,width=4.5cm,
                      ytick={0,0.2,0.4,0.6,0.8,1},
                      ymax=1.1]
\nextgroupplot[title={$n=200$}]\addplot [domain=0:.5, color=black] coordinates { 
(0,0.04)
(0.1,0.052)
(0.2,0.114)
(0.3,0.18)
(0.4,0.314)
(0.5,0.472)
}; 
\addplot [domain=0:.5, color=brightgreen] coordinates { 
(0,0.052)
(0.1,0.056)
(0.2,0.138)
(0.3,0.256)
(0.4,0.43)
(0.5,0.666)
}; 
\addplot [domain=0:.5, color=blue] coordinates { 
(0,0.054)
(0.1,0.04)
(0.2,0.114)
(0.3,0.222)
(0.4,0.344)
(0.5,0.514)
}; 
\addplot [domain=0:.5, color=orange] coordinates { 
(0,0.016)
(0.1,0.016)
(0.2,0.008)
(0.3,0.026)
(0.4,0.032)
(0.5,0.032)
}; 
\addplot [domain=0:.5, color=red] coordinates { 
(0,0.032)
(0.1,0.058)
(0.2,0.116)
(0.3,0.216)
(0.4,0.374)
(0.5,0.552)
}; 
\addplot [domain=0:.5, color=grey] coordinates { 
(0,0.05)
(0.1,0.05)
(0.2,0.05)
(0.3,0.05)
(0.4,0.05)
(0.5,0.05)
}; 
\nextgroupplot[title={$n=400$}]\addplot [domain=0:.5, color=black] coordinates { 
(0,0.05)
(0.1,0.0880000000000001)
(0.2,0.176)
(0.3,0.39)
(0.4,0.664)
(0.5,0.888000000000001)
}; 
\addplot [domain=0:.5, color=brightgreen] coordinates { 
(0,0.058)
(0.1,0.102)
(0.2,0.226)
(0.3,0.494)
(0.4,0.796000000000001)
(0.5,0.940000000000001)
}; 
\addplot [domain=0:.5, color=blue] coordinates { 
(0,0.044)
(0.1,0.0780000000000001)
(0.2,0.19)
(0.3,0.42)
(0.4,0.692000000000001)
(0.5,0.888000000000001)
}; 
\addplot [domain=0:.5, color=orange] coordinates { 
(0,0.01)
(0.1,0.014)
(0.2,0.02)
(0.3,0.028)
(0.4,0.02)
(0.5,0.038)
}; 
\addplot [domain=0:.5, color=red] coordinates { 
(0,0.048)
(0.1,0.0960000000000001)
(0.2,0.198)
(0.3,0.476)
(0.4,0.772000000000001)
(0.5,0.944000000000001)
}; 
\addplot [domain=0:.5, color=grey] coordinates { 
(0,0.05)
(0.1,0.05)
(0.2,0.05)
(0.3,0.05)
(0.4,0.05)
(0.5,0.05)
}; 
\nextgroupplot[title={$n=800$}]\addplot [domain=0:.5, color=black] coordinates { 
(0,0.038)
(0.05,0.058)
(0.1,0.132)
(0.15,0.258)
(0.2,0.386)
(0.25,0.594)
(0.3, 0.76)
(0.4, 0.968)
(0.5, 0.996)
}; 
\label{plots:plot1VDW5D8} 
\addplot [domain=0:.5, color=brightgreen] coordinates { 
(0,0.034)
(0.05,0.058)
(0.1,0.162)
(0.15,0.33)
(0.2,0.504)
(0.25,0.708)
(0.3, 0.882)
(0.4, 0.984)
(0.5, 1)
}; 
\label{plots:plot2VDW5D8} 
\addplot [domain=0:.5, color=blue] coordinates { 
(0,0.062)
(0.05,0.056)
(0.1,0.162)
(0.15,0.278)
(0.2,0.466)
(0.25,0.628)
(0.3, 0.848)
(0.4, 0.982)
(0.5, 1)
}; 
\label{plots:plot3VDW5D8} 
\addplot [domain=0:.5, color=red] coordinates { 
(0,0.042)
(0.05,0.052)
(0.1,0.152)
(0.15,0.272)
(0.2,0.474)
(0.25,0.718)
(0.3, 0.854)
(0.4, 0.984)
(0.5, 1)
}; 
\label{plots:plot4VDW5D8} 
\addplot [domain=0:.5, color=orange] coordinates { 
(0,0.016)
(0.05,0.006)
(0.1,0.014)
(0.15,0.01)
(0.2,0.02)
(0.25,0.024)
(0.3,0.018)
(0.4,0.034)
(0.5,0.04)

}; 
\label{plots:plot5VDW5D8} 
\addplot [domain=0:.5, color=grey] coordinates { 
(0,0.05)
(0.05,0.05)
(0.1,0.05)
(0.15,0.05)
(0.2,0.05)
(0.25,0.05)
(0.5,0.05)
}; 
\label{plots:plot6VDW5D8} 

            \end{groupplot}      
\path (myplot c1r1.south west|-current bounding box.south)--
      coordinate(legendpos)
      (myplot c3r1.south east|-current bounding box.south);
\matrix[
    matrix of nodes,
    anchor=south,
    draw,
    inner sep=0.2em,
    draw
  ]at([yshift=-6ex]legendpos)
  {
     \ref{plots:plot3VDW}& $\left(\tenq{T}^{(n)}_{\text{\rm vdW}{\pm}}\right)^2$&[5pt]
        \ref{plots:plot4VDW}&$\left(\tenq{T}^{(n)}_{\text{\rm vdW}{\scriptscriptstyle\square}}\right)^2$ &[5pt]
    \ref{plots:plot2VDW}& $\left(\tenq{T}^{(n)}_{\text{\rm vdW}{\pm}{\mathcal N}}\right)^2$&[5pt]
             \ref{plots:plot1VDW}& $\left(\tenq{T}^{(n)}_{\text{\rm vdW}{\scriptscriptstyle\square}{\mathcal N}}\right)^2$&[5pt]
              \ref{plots:plotHVDW}& $T^2$
 \\};
\end{tikzpicture}
\caption{\small\slshape Rejection frequencies,   for   5-dimensional   samples with independent Cauchy marginals,  (see~Section~\ref{sec524}) and various sample sizes, of Hotelling's test based on~$T^2$ and the van der Waerden tests based on $\left(\tenq{T}^{(n)}_{\text{\rm vdW}{\pm}}\right)^2\vspace{-1mm}$,  $\left(\tenq{T}^{(n)}_{\text{\rm vdW}{\scriptscriptstyle\square}}\right)^2$, $\left(\tenq{T}^{(n)}_{\text{\rm vdW}{\pm}{\mathcal N}}\right)^2$, 
and~$\left(\tenq{T}^{(n)}_{\text{\rm vdW}{\scriptscriptstyle\square}{\mathcal N}}\right)^2\vspace{1.5mm}$, respectively,  as  functions of the shift $\eta$;   $N=500$ replications.} 
\label{fig: CompCauchyVDW5}
\end{figure}



\subsubsection{Samples with independent Cauchy marginals}\label{sec524}\vspace{-3mm}

Same  Cauchy samples   as in Section~\ref{sec424}. 
Rejection frequencies over~$N=~\!500$ 
 replications are
 shown (as  functions of $\eta$) in  Figure~\ref{fig: CompCauchyVDW5}. The conclusions drawn for~$d=2$ still hold, with a very slight superiority of the ``direct transportation'' test~$\left(\tenq{T}^{(n)}_{\text{\rm vdW}{\pm}{\mathcal N}}\right)^2$ over  the Gaussian score ones $\left(\tenq{T}^{(n)}_{\text{\rm vdW}{\pm}}\right)^2$ and $\left(\tenq{T}^{(n)}_{\text{\rm vdW}{\scriptscriptstyle\square}}\right)^2$.

\section{Conclusions}
While confirming the advantages and excellent performance of rank tests over their daily practice pseudo-Gaussian counterparts, the simulations of the previous sections provide  empirical answers to several  questions of great practical importance. 

The  choice of the grid (whether spherical ($\mathfrak G$i), cubic ($\mathfrak G$ii), or Gaussian (($\mathfrak G$iii) and~($\mathfrak G$iv)) seems to have relatively  little impact on the performance of the corresponding Wilcoxon tests (the only significant case being the Cauchy one), and no impact at all on the performance of van der Waerden tests. In particular, there is no evidence that Wilcoxon tests based on spherical grids~($\mathfrak G$i) are preferable under spherical distributions while Wilcoxon tests based on  cubic grids ($\mathfrak G$ii)  are preferable under distributions with  independent components: see, e.g., the Cauchy case   (Sections~\ref{sec416} and~\ref{sec414}).


\bibliographystyle{apa-good}

\bibliography{Location_MvRanks}

\end{document}